\documentclass[english]{article}
\usepackage[T1]{fontenc}
\usepackage[utf8]{inputenc}
\usepackage{geometry}
\usepackage{amsmath}
\usepackage{amssymb}
\usepackage{stackrel}
\usepackage{graphicx}

\makeatletter

\providecommand{\tabularnewline}{\\}

\usepackage{amsfonts}
\usepackage{amsthm}
\usepackage{bbold}
\usepackage[all]{xy}
\usepackage{babel}

\makeatother

\usepackage{babel}
\begin{document}
\begin{center}
\bigskip{}
\bigskip{}
\bigskip{}
\bigskip{}
\bigskip{}
\bigskip{}
\textbf{\textsc{\LARGE{}a ``magic''\ approach to octonionic }}\\
\textbf{\textsc{\LARGE{}Rosenfeld spaces }}{\LARGE\par}
\par\end{center}

\bigskip{}

\begin{center}
{\large{}12 December 2022}{\large\par}
\par\end{center}

\bigskip{}

\begin{center}
{\large{}Alessio Marrani$^{a}$, Daniele Corradetti$^{b}$, David
Chester$^{c}$, }{\large\par}
\par\end{center}

\begin{center}
{\large{}Raymond Aschheim$^{c}$, Klee Irwin$^{c}$}{\large\par}
\par\end{center}

\bigskip{}

\begin{center}
\emph{$a$~Instituto de Física Teorica, Departamento de Física, }\\
\emph{Universidad de Murcia, Campus de Espinardo, E-30100, Spain}\\
\emph{email: }\texttt{alessio.marrani@um.es} 
\par\end{center}

\begin{center}
\emph{$b$ Universidade do Algarve, Departamento de Matemática, }\\
\emph{Campus de Gambelas, 8005-139 Faro, Portugal}\\
\emph{email: }\texttt{a55944@ualg.pt}
\par\end{center}

\begin{center}
\emph{$c$ Quantum Gravity Research, }\\
\emph{Topanga Canyon Rd 101 S., California CA 90290, USA}\\
\emph{email: }\texttt{DavidC@QuantumGravityResearch.org}; \\
\texttt{Raymond@QuantumGravityResearch.org}; \\
\texttt{Klee@QuantumGravityResearch.org}
\par\end{center}

\bigskip{}
\bigskip{}
\bigskip{}
\bigskip{}

In his study on the geometry of Lie groups, Rosenfeld postulated a
strict relation between all real forms of exceptional Lie groups and
the isometries of projective and hyperbolic spaces over the (rank-2)
tensor product of Hurwitz algebras taken with appropriate conjugations.
Unfortunately, the procedure carried out by Rosenfeld was not rigorous,
since many of the theorems he had been using do not actually hold
true in the case of algebras that are not alternative nor power-associative.
A more rigorous approach to the definition of all the planes presented
more than thirty years ago by Rosenfeld in terms of their isometry
group, can be considered within the theory of coset manifolds, which
we exploit in this work, by making use of all real forms of Magic
Squares of order three and two over Hurwitz normed division algebras
and their split versions. Within our analysis, we find 7 pseudo-Riemannian
symmetric coset manifolds which seemingly cannot have any interpretation
within Rosenfeld's framework. We carry out a similar analysis for
Rosenfeld lines, obtaining that there are a number of pseudo-Riemannian
symmetric cosets which do not have any interpretation \textit{à la
Rosenfeld}.

\newpage{}

\tableofcontents{}

\newpage{}

\section{Introduction}

Around the half of the XX century, geometric investigations on the
octonionic plane gave rise to a fruitful mathematical activity, culminating
into the formulation of \textit{Tits-Freudenthal Magic Square} \cite{Freud-1965,Tits}.
The Magic Square is an array of Lie algebras, whose entries are obtained
from two Hurwitz algebras\footnote{Even though the original construction by Tits was not symmetric in
the two entries (cf. the third of (\ref{msconstructions})), the square
turned out to be symmetric and was, therefore, dubbed as ``magic''.} and enjoy multiple geometric and algebraic interpretations \cite{pre-BS,BS,Santander,Santander-2,EldMS2,EldMS1},
as well as physical applications \cite{GST-1,Wissanji,Squaring-Magic}.

The entries $\mathfrak{m}_{3}(\mathbb{A}_{1},\mathbb{A}_{2})$ of
the Tits-Freudenthal Magic Square can be defined equivalently \cite{BS}
using the Tits \cite{Tits}, Barton-Sudbery \cite{BS} (also cf. \cite{Evans:2009ed}),
and Vinberg \cite{Vinberg} constructions\footnote{The ternary algebra approach of \cite{Kantor} was generalised by
Bars and Günaydin in \cite{Bars-Gun} to include super Lie algebras.
Generalizations to affine, hyperbolic and further extensions of Lie
algebras have been considered in \cite{Palmkvist}.}, 
\begin{equation}
\mathfrak{m}_{3}(\mathbb{A}_{1},\mathbb{A}_{2})=\left\{ \begin{array}{l}
\mathfrak{der}(\mathbb{A}_{1})\oplus\mathfrak{der}(\mathfrak{J}_{3}(\mathbb{A}_{2}))\,\dot{+}\,\text{Im}\mathbb{A}_{1}\otimes\mathfrak{J}_{3}^{\prime}(\mathbb{A}_{2});\\[5pt]
~\\
\mathfrak{der}(\mathbb{A}_{1})\oplus\mathfrak{der}(\mathbb{A}_{2})\,\dot{+}\,\mathfrak{sa}(3,\mathbb{A}_{1}\otimes\mathbb{A}_{2});\\
~\\
\mathfrak{tri}(\mathbb{A}_{1})\oplus\mathfrak{tri}(\mathbb{A}_{2})\,\dot{+}\,3(\mathbb{A}_{1}\otimes\mathbb{A}_{2}),
\end{array}\right.\label{msconstructions}
\end{equation}
where $\mathbb{A}_{1}$ and $\mathbb{A}_{2}$ are normed Hurwitz division
algebras $\mathbb{R},\mathbb{C},\mathbb{H},\mathbb{O}$ or the split
counterparts $\mathbb{C}_{s},\mathbb{H}_{s},\mathbb{O}_{s}$. Here,
$\mathfrak{J}_{3}(\mathbb{A})$ denotes the Jordan algebra of $3\times3$
Hermitian matrices over $\mathbb{A}$ and $\mathfrak{J}_{3}^{\prime}(\mathbb{A})$
its subspace of traceless elements. The space of anti-Hermitian traceless
$n\times n$ matrices over $\mathbb{A}_{1}\otimes\mathbb{A}_{2}$
is denoted $\mathfrak{sa}(n,\mathbb{A}_{1}\otimes\mathbb{A}_{2})$.
Details of the commutators and the isomorphisms between these Lie
algebras can be found e.g. in \cite{BS}. The choice of division or
split $\mathbb{A}_{1},\mathbb{A}_{2}$ yields different real forms
: for the Tits construction, further possibilities are given by allowing
the Jordan algebra to be Lorentzian, denoted $\mathfrak{J}_{2,1}$,
as described in \cite{Squaring-Magic} (see Sec. \ref{Eucl-Lor} below);
equivalently, one can introduce overall signs in the definition of
the commutators between distinct components in the Vinberg or Barton-Sudbery
constructions, as described in \cite{Magic-Pyramids}. For a complete
listing of all possibilities\footnote{Concerning the extension of Freudenthal-Tits Magic Square to various
types of ``intermediate''\ algebras, the extension to the sextonions
$\mathbb{S}$ (between $\mathbb{H}$ and $\mathbb{O}$) is due to
Westbury \cite{Westbury}, and later developed by Landsberg and Manivel
\cite{Lands-Maniv,Lands-Maniv-sext}, whereas further extension to
the tritonions $\mathbb{T}$ (between $\mathbb{C}$ and $\mathbb{H}$)
was discussed by Borsten and one of the present authors in \cite{Kind-of-Magic}.}, see \cite{Squaring-Magic}. Recently \cite{WDM}, Wilson, Dray and
Manogue gave a new construction of the Lie algebra $\mathfrak{e}_{8}$
in terms of $3\times3$ matrices such that the Lie bracket has a natural
description as matrix commutator; this led to a new interpretation
of the Freudenthal-Tits Magic Square of Lie algebras\footnote{Magic Squares \textit{of order 2} of Lie groups have been investigated
in \cite{Dray-Huerta-Kincaid}. Interestingly, the Freudenthal-Tits
Magic Square of Lie groups was observed to be non-symmetric since
\cite{Yokota-2}.}, acting on themselves by commutation.

While Tits was more interested in algebraic aspects, Freudenthal interpreted
every row of the Magic Square in terms of a different kind of geometry,
i.e. elliptic, projective, symplectic and meta-symplectic for the
first, second, third and fourth row, respectively \cite{Lands-Maniv}.
This means that, while the algebras taken in account were always the
Hurwitz algebras $\mathbb{R},\mathbb{C},\mathbb{H}$ and $\mathbb{O}$,
all Lie groups were arising from considering different type of transformations
such as isometries, collineations and homographies of the plane \cite{Freud-1965}.
On the other hand, Rosenfeld conceived every entry of the $4\times4$
array of the Magic Square as the Lie algebra of the (global) isometry
group of a ``generalized''\ projective plane, later called \textit{Rosenfeld
plane} \cite{Rosenf98,RosenfeldGroup2-1}. Thus, while in Freudenthal's
framework the projective plane was always the same but, depending
on the considered row of the Magic Square, the kind of transformations
was changed; in Rosenfeld's picture the group of transformations was
always the isometry group, and different Lie groups appeared considering
different planes over tensor product of Hurwitz algebras. 

To be more precise, since the exceptional Lie group $\text{F}_{4}$
is the isometry group of the octonionic projective plane $\mathbb{O}P^{2}$,
then Rosenfeld regarded $\text{E}_{6}$ as the isometry group of the
\textit{bioctonionic ``projective''\ plane} $\mathbb{\left(C\otimes\mathbb{O}\right)}P^{2}$,
$\text{E}_{7}$ as the isometry group of the \textit{quaternoctonionic
``projective''\ plane} $\mathbb{\left(H\otimes\mathbb{O}\right)}P^{2}$,
and $\text{E}_{8}$ as the isometry group of the \textit{octooctonionic
``projective''\ plane} $\mathbb{\left(O\otimes\mathbb{O}\right)}P^{2}$.
Despite Rosenfeld's suggestive interpretation, it was soon realized
that the planes identified by Rosenfeld did not satisfy projective
axioms, and this explains the quotes in the term ``projective'',
which would have to be intended only in a vague sense. It was the
work of Atsuyama, followed by Landsberg and Manivel, to give a rigorous
description of the geometry arising from Rosenfeld's approach \cite{Atsuyama,Lands-Maniv,Atiyah-Bernd}.
Indeed (unlike $\mathbb{O\simeq R\otimes O}$) $\mathbb{C\otimes O}$,
$\mathbb{H\otimes O}$ and $\mathbb{O\otimes O}$ are not division
algebras, thus preventing a direct projective construction; moreover
(unlike $\mathbb{C\otimes O}$), Hermitian $3\times3$ matrices over
$\mathbb{H\otimes O}$ or $\mathbb{O\otimes O}$ do not form a simple
Jordan algebra, so the usual identification of points (lines) with
trace 1 (2) projection operators cannot be made \cite{Baez}. Nonetheless,
they are in fact geometric spaces, generalising projective spaces,
known as ``buildings'', on which (the various real forms of) exceptional
Lie groups act as isometries. Buildings where originally introduced
by Tits in order to provide a geometric approach to simple Lie groups,
in particular the exceptional cases, but have since had far reaching
implications; see, for instance, \cite{Build-1,Build-2} and Refs.
therein.

More specifically, Rosenfeld noticed that, over $\mathbb{R}$, the
isometry group of the octonionic projective plane was the compact
form of $\text{F}_{4}$, i.e. $\text{F}_{4\left(-52\right)}$, while
the other, non-compact real forms, i.e. the split form $\text{F}_{4\left(4\right)}$
and $\text{F}_{4\left(-20\right)}$, were obtained as isometry groups
of the \textit{split-octonionic projective plane} $\mathbb{O}_{s}P^{2}$
and of the \textit{octonionic hyperbolic plane} $\mathbb{O}H^{2}$
\cite{Corr-Notes-Octo,RealF}, respectively. Thence, he proceeded
in relating all real forms of exceptional Lie groups with projective
and hyperbolic planes over tensorial products of Hurwitz algebras
\cite{RosenfeldGroup2-1,Rosenf98}.

Despite being very insightful, to the best of our knowledge Rosenfeld's
approach was never formulated in a systematic way, and a large part
of the very Rosenfeld planes were never really rigorously defined.
In the present investigation, we give a systematic and explicit definition
of all octonionic\footnote{The quaternionic, complex and real \textit{``Rosenfeld\
planes''} and \textit{``Rosenfeld lines''} can be obtained as proper
sub-manifolds of their octonionic counterparts, so we will here focus
only on the latter ones.} ``\textit{Rosenfeld planes}''\ as coset manifolds, by exploiting
all real forms of Magic Squares of order 3 over Hurwitz algebras and
their split versions, thus all possible real forms of Freudenthal-Tits
Magic Square. In doing so, we will also find a total of 10 coset manifolds,
namely 7 ``planes''\ and 3 ``lines'', whose definition is straightforward
according to our procedure, but that apparently do not have any interpretation
according to Rosenfeld's approach.

The plan of the paper is as follows. In Sec. \ref{1} we will resume
all the algebraic machinery used in this paper, i.e. Hurwitz algebras,
their triality and derivation symmetries, Euclidean and Lorentzian
cubic Jordan algebras, Magic Squares of order 3 and 2, etc. In Sec.
\ref{2} we will then present the ``magic''\ formulæ\ which we
will use in order to define the octonionic Rosenfeld planes in Sec.
\ref{3}; then, we will do the same for octonionic Rosenfled lines
in Secs. \ref{pre-4} and \ref{4}. Some final remarks and an outlook
are given in Sec. \ref{5}.

\section{\label{1}Real forms of Magic Squares of order 3 and 2}

\subsection{Hurwitz algebras and their split versions}

\begin{figure}
\centering{}\includegraphics[scale=0.06]{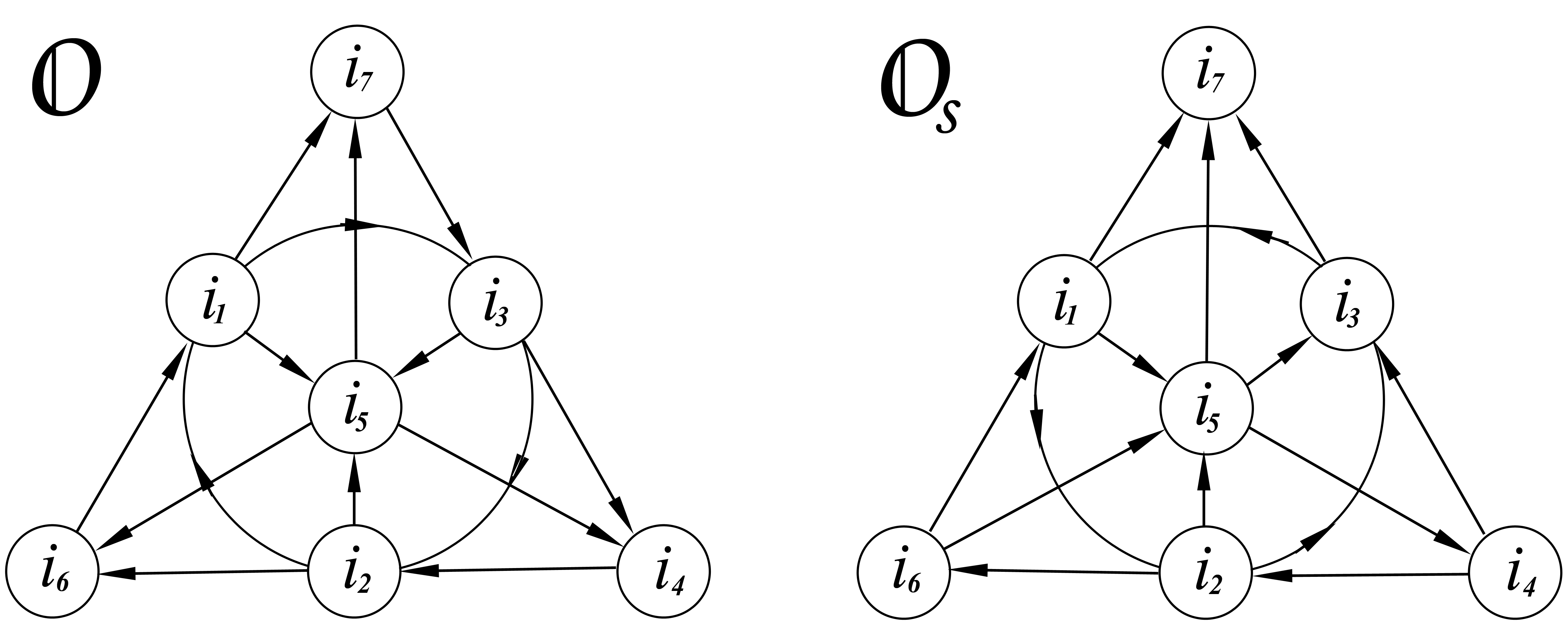}\caption{\label{fig:octonion fano plane} Multiplication rule of octonions
$\mathbb{O}$ (\emph{left}) and of split-octonions $\mathbb{O}_{s}$
(\emph{right}) as real vector space $\mathbb{R}^{8}$ in the basis
$\left\{ \text{i}_{0}=1,\text{i}_{1},...,\text{i}_{7}\right\} $.
In the case of the octonions $\text{i}_{0}^{2}=1$ and $\text{i}_{k}^{2}=-1$
for $k=1...7$, while in the case of split-octonions $\text{i}_{k}^{2}=+1$,
for $k=1,2,3$ and $\text{i}_{k}^{2}=-1$ for $k\protect\neq1,2,3$.}
\end{figure}

Let the \textit{octonions} $\mathbb{O}$ be the only unital, non-associative,
normed division algebra with $\mathbb{R}^{8}$ decomposition of $x\in\mathbb{O}$
given by 
\begin{equation}
x=\sum_{k=0}^{7}x_{k}\text{i}_{k},\label{eq:decomposizione}
\end{equation}
where $\left\{ \text{i}_{0}=1,\text{i}_{1},...,\text{i}_{7}\right\} $
is a basis of $\mathbb{R}^{8}$ and the multiplication rules are mnemonically
encoded in the \textit{Fano plane} (see left side of Fig. \ref{fig:octonion fano plane}),
along with $\text{i}_{k}^{2}=-1$ for $k=1,...,7$. Let the \textit{norm}
$N\left(x\right)$ be defined as 
\begin{equation}
N\left(x\right):=\sum_{k=0}^{7}\left(x_{k}\text{i}_{k}\right)^{2},\label{norm}
\end{equation}
and its polarisation as $\left\langle x,y\right\rangle :=N\left(x+y\right)-N\left(x\right)-N\left(y\right).$
Then the \textit{real part} of $x$ is $\mathfrak{R}\left(x\right):=\left\langle x,1\right\rangle $
and $\overline{x}:=2\left\langle x,1\right\rangle -x.$ We then have
that $N\left(x\right)=\overline{x}x,\,\,\,\left\langle x,y\right\rangle =\overline{x}y+\overline{y}x,$and
that 
\begin{equation}
N\left(xy\right)=N\left(x\right)N\left(y\right),
\end{equation}
or, in other words, that octonions are a \textit{composition algebra}
with respect to the Norm $N$ defined by (\ref{norm}).

We now define as the algebra of \textit{quaternions} $\mathbb{H}$
as the subalgebra of $\mathbb{O}$ generated by the elements $\left\{ \text{i}_{0}=1,\text{i}_{1},\text{i}_{2},\text{i}_{3}\right\} $
and the algebra of \textit{complex numbers} $\mathbb{C}$ as the subalgebra
of $\mathbb{O}$ generated by the elements $\left\{ \text{i}_{0}=1,\text{i}_{1}\right\} $.

Thus, the \textit{split-octonions} $\mathbb{O}_{s}$ can be defined
as the only unital, non-associative, normed algebra with $\mathbb{R}^{8}$
decomposition of $x\in\mathbb{O}_{s}$ given by (\ref{eq:decomposizione})
where, again, $\left\{ \text{i}_{0}=1,\text{i}_{1},...,\text{i}_{7}\right\} $
is a basis of $\mathbb{R}^{8}$, but the multiplication rules are
encoded in a different variant of the Fano plane, given by the right
side of Fig. \ref{fig:octonion fano plane}, along with $\text{i}_{k}^{2}=1$
for $k=1,...,3$ and $\text{i}_{k}^{2}=-1$ for $k=4,5,6,7$ . Norm
and conjugation are defined as in the octonionic case and, as for
the division octonions, we obtain that split octonions are a \textit{composition
algebra} with respect to the Norm $N$ defined by (\ref{norm}). We
then define the algebra of \textit{split-quaternions} $\mathbb{H}_{s}$
as the subalgebra of $\mathbb{O}_{s}$ generated by $\left\{ \text{i}_{0}=1,\text{i}_{1},\text{i}_{4},\text{i}_{5}\right\} $
and the \textit{split-complex} (or \textit{hypercomplex}) \textit{numbers}
$\mathbb{C}_{s}$ as the subalgebra of $\mathbb{O}_{s}$ generated
by $\left\{ \text{i}_{0}=1,\text{i}_{4}\right\} $.

\subsection{Trialities and derivations}

Let $\mathbb{A}$ be an algebra and $\text{End}\left(\mathbb{A}\right)$
the associative algebra of its linear endomorphisms. Then, the triality
algebra $\mathfrak{tri}\left(\mathbb{A}\right)$ is the Lie subalgebra
of $\mathfrak{so}\left(\mathbb{A}\right)\oplus\mathfrak{so}\left(\mathbb{A}\right)\oplus\mathfrak{so}\left(\mathbb{A}\right)$
defined as 
\begin{equation}
\mathfrak{tri}\left(\mathbb{A}\right):=\left\{ \left(A,B,C\right)\in\text{End}\left(\mathbb{A}\right)^{\oplus3}:A\left(xy\right)=B\left(x\right)y+xC\left(y\right),\forall x,y\in\mathbb{A}\right\} ,
\end{equation}
where the Lie bracket are those inherited as a subalgebra. Derivations
are a special case of trialities of the form $\left(A,A,A\right)$
and therefore they again form a Lie algebra, defined as 
\begin{equation}
\mathfrak{der}\left(\mathbb{A}\right):=\left\{ A\in\text{End}\left(\mathbb{A}\right):A\left(xy\right)=A\left(x\right)y+xA\left(y\right),\forall x,y\in\mathbb{A}\right\} ,
\end{equation}
where the bracket is given by the commutator. In the case of $\mathbb{A}=\mathbb{C},\mathbb{H},\mathbb{O}$
we intended $\mathfrak{so}\left(\mathbb{A}\right)$ as $\mathfrak{so}\left(2\right)$,
$\mathfrak{so}\left(4\right)$ and $\mathfrak{so}\left(8\right)$
respectively; while, for their split companions, we intended $\mathfrak{so}\left(\mathbb{A}_{s}\right)$
as $\mathfrak{so}\left(1,1\right)$, $\mathfrak{so}\left(2,2\right)$
and $\mathfrak{so}\left(4,4\right)$. Triality and derivation Lie
algebras over Hurwitz algebras are summarized in Table \ref{tab:Triality-algebras}.
\begin{table}
\begin{centering}
\begin{tabular}{|c|c|c|c|c|}
\hline 
 & $\mathbb{R}$ & $\mathbb{C}$ & $\mathbb{H}$ & $\mathbb{O}$\tabularnewline
\hline 
\hline 
$\mathbb{\mathfrak{tri}\left(\mathbb{A}\right)}$ & $\textrm{Ø}$ & $\mathfrak{u}\left(1\right)\oplus\mathfrak{u}\left(1\right)$ & $\mathfrak{so}\left(3\right)\oplus\mathfrak{so}\left(3\right)\oplus\mathfrak{so}\left(3\right)$ & $\mathfrak{so}\left(8\right)$\tabularnewline
\hline 
$\mathbb{\mathfrak{tri}}\left(\mathbb{A}_{s}\right)$ & $\textrm{Ø}$ & $\mathfrak{so}\left(1,1\right)\oplus\mathfrak{so}\left(1,1\right)$ & $\mathfrak{sl}_{2}\left(\mathbb{R}\right)\oplus\mathfrak{sl}_{2}\left(\mathbb{R}\right)\oplus\mathfrak{sl}_{2}\left(\mathbb{R}\right)$ & $\mathfrak{so}\left(4,4\right)$\tabularnewline
\hline 
$\mathfrak{der}\left(\mathbb{A}\right)$ & $\textrm{Ø}$ & $\textrm{Ø}$ & $\mathfrak{so}\left(3\right)$ & $\mathfrak{g}_{2\left(-14\right)}$\tabularnewline
\hline 
$\mathfrak{der}\left(\mathbb{A}_{s}\right)$ & $\textrm{Ø}$ & $\textrm{Ø}$ & $\mathfrak{sl}_{2}\left(\mathbb{R}\right)$ & $\mathfrak{g}_{2\left(2\right)}$\tabularnewline
\hline 
$\mathfrak{A}\left(\mathbb{A}\right)$ & $\textrm{Ø}$ & $\mathfrak{u}_{1}$ & $\mathfrak{su}_{2}$ & $\textrm{Ø}$\tabularnewline
\hline 
$\mathfrak{A}\left(\mathbb{A}_{s}\right)$ & $\textrm{Ø}$ & $\mathfrak{so}_{1,1}$ & $\mathfrak{sl}_{2}\left(\mathbb{R}\right)$ & $\textrm{Ø}$\tabularnewline
\hline 
\end{tabular}
\par\end{centering}
\caption{\label{tab:Triality-algebras}Triality and derivation algebras of
Hurwitz algebras. Moreover we also added the algebra $\mathfrak{A}\left(\mathbb{A}\right):=\mathfrak{tri}\left(\mathbb{A}\right)\ominus\mathfrak{so}\left(\mathbb{A}\right),$
that will be used in Sec. 3.2 in the definition of the Rosenfeld planes.}
\end{table}

\subsection{\label{Eucl-Lor}Euclidean and Lorentzian simple rank-3 Jordan algebras}

Let $\mathbb{A}$ be a composition algebra and $H_{n}\left(\mathbb{A}\right)$
be the set of Hermitian $n\times n$ matrices with elements in $\mathbb{A}$
such that
\begin{equation}
X^{\dagger}:=\overline{X}^{t}=X,\label{Herm}
\end{equation}
where the conjugation is the one pertaining to $\mathbb{A}$ itself.
When\footnote{Excluding $\mathbb{A}=\mathbb{O}$ and $\mathbb{O}_{s}$, one can
actually consider\textit{\ any} $n\in\mathbb{N}$.} $n=2$ or $3$, we define the simple, rank-$n$ \textit{Euclidean
Jordan algebra }$\mathfrak{J}_{n}\left(\mathbb{A}\right)$ as the
commutative algebra over $H_{n}\left(\mathbb{A}\right)$ with the
\textit{Jordan product} $\circ$ defined by the anticommutator 
\begin{equation}
X\circ Y:=\frac{1}{2}\left(XY+YX\right)=:\left\{ X,Y\right\} ,~\forall X,Y\in H_{n}\left(\mathbb{A}\right),
\end{equation}
where juxtaposition denotes the standard 'rows-by-columns' matrix
product. Moreover, we define $\mathfrak{J}_{n}^{\prime}\left(\mathbb{A}\right)$
as the subspace of the Jordan algebra orthogonal to the identity $I$
(namely, the subspace of traceless matrices in $\mathfrak{J}_{n}\left(\mathbb{A}\right)$),
endowed with the product inherited\footnote{Note that $\mathfrak{J}_{n}^{\prime}\left(\mathbb{A}\right)$ is not
closed under the restriction of the \textit{Jordan product} $\circ$
of $\mathfrak{J}_{n}\left(\mathbb{A}\right)$ to $\mathfrak{J}_{n}^{\prime}\left(\mathbb{A}\right)$
itself. However, one can deform $\circ$ into the \textit{Michel-Radicati
product} $\ast$, defined as $X\ast Y:=X\circ Y-\frac{\text{Tr}\left(X\circ Y\right)}{n}I$
$\forall X,Y\in\mathfrak{J}_{n}^{\prime}\left(\mathbb{A}\right)$,
under which $\mathfrak{J}_{n}^{\prime}\left(\mathbb{A}\right)$ is
closed \cite{MR70,MR73}.} from $\mathfrak{J}_{n}\left(\mathbb{A}\right)$ (see e.g. \cite[p.8]{BS}).
In particular, for $n=3$ the simple, rank-$3$ \textit{Euclidean
Jordan algebra }$\mathfrak{J}_{3}\left(\mathbb{A}\right)$ has elements
\begin{equation}
X:=\left(\begin{array}{ccc}
a & x_{1} & x_{2}\\
\overline{x}_{1} & b & x_{3}\\
\overline{x}_{2} & \overline{x}_{3} & c
\end{array}\right)\in\mathfrak{J}_{3}\left(\mathbb{A}\right),
\end{equation}
where $x_{1},x_{2},x_{3}\in\mathbb{A}$ and $a,b,c\in\mathbb{R}$.

The Hermiticity condition (\ref{Herm}) can be generalized by inserting
a pseudo-Euclidean metric :
\begin{eqnarray}
 &  & \eta X^{\dagger}\eta=X,\label{Herm-2}\\
 &  & \eta:=\text{diag}\left(\underset{p~\text{times}}{\underbrace{-1,...,-1}},\underset{n-p~\text{times}}{\underbrace{1,...,1}}\right).
\end{eqnarray}
Correspondingly, one defines the simple, rank-$n$\textit{ pseudo-Euclidean
Jordan algebra }$\mathfrak{J}_{n-p,p}\left(\mathbb{A}\right)\simeq\mathfrak{J}_{p,n-p}\left(\mathbb{A}\right)$.
In particular, by setting $n=3$, the simple, rank-$3$ \textit{Lorentzian
Jordan algebra }$\mathfrak{J}_{2,1}\left(\mathbb{A}\right)\simeq\mathfrak{J}_{1,2}\left(\mathbb{A}\right)$,
whose elements read as 
\begin{equation}
X:=\left(\begin{array}{ccc}
a & x_{1} & x_{2}\\
-\overline{x}_{1} & b & x_{3}\\
-\overline{x}_{2} & \overline{x}_{3} & c
\end{array}\right)\in\mathfrak{J}_{1,2}\left(\mathbb{A}\right),
\end{equation}
where $x_{1},x_{2},x_{3}\in\mathbb{A}$ and $a,b,c\in\mathbb{R}$.

\subsection{Real forms of the Freudenthal-Tits Magic Square}

As anticipated in (\ref{msconstructions}), given $\mathbb{A}_{1}$
and $\mathbb{A}_{2}$ two composition algebras, the corresponding
entry of Tits-Freudenthal Magic Square $\mathfrak{m}_{3}\left(\mathbb{A}_{1},\mathbb{A}_{2}\right)$
is a Lie algebra that can be realized in \textit{at least} three different,
but equivalent ways, on which we will now briefly comment (without
the explicit definition of the Lie brackets, that can be found e.g.
in \cite{BS}).
\begin{enumerate}
\item Tits construction \cite{Tits} is given by the first line of the r.h.s.
of (\ref{msconstructions}) :
\begin{equation}
\mathfrak{m}_{3}\left(\mathbb{A}_{1},\mathbb{A}_{2}\right):=\mathfrak{der}(\mathbb{A}_{1})\oplus\mathfrak{der}(\mathfrak{J}_{3}(\mathbb{A}_{2}))\,\dot{+}\,\text{Im}\mathbb{A}_{1}\otimes\mathfrak{J}_{3}^{\prime}(\mathbb{A}_{2}),\label{Tits'}
\end{equation}
thus involving the derivation Lie algebra of $\mathbb{A}_{1}$ and
the derivation Lie algebra of the simple, rank-$3$ Euclidean Jordan
algebra $\mathfrak{J}_{3}\left(\mathbb{A}_{2}\right)$ over $\mathbb{A}_{2}$.
Tits construction, despite being not manifestly symmetric under the
exchange $\mathbb{A}_{1}\leftrightarrow\mathbb{A}_{2}$, is the most
general of the three constructions presented here, since it holds
for any alternative algebras $\mathbb{A}_{1}$ and $\mathbb{A}_{2}$,
as long as it is possible to define a Jordan algebra over $\mathbb{A}_{2}$
itself. As mentioned above, Tits construction can be generalized to
$\mathfrak{J}_{2,1}(\mathbb{A}_{2})$, thus yielding\footnote{Since $\mathfrak{J}_{2,1}(\mathbb{A}_{2})\simeq\mathfrak{J}_{1,2}(\mathbb{A}_{2})$,
it holds that $\mathfrak{m}_{2,1}(\mathbb{A}_{1},\mathbb{A}_{2})\simeq\mathfrak{m}_{1,2}(\mathbb{A}_{1},\mathbb{A}_{2})$.} $\mathfrak{m}_{1,2}\left(\mathbb{A}_{1},\mathbb{A}_{2}\right)$.
\item A more symmetric approach was pursued by Vinberg \cite{Vinberg},
who obtained the formula given by the second line of the r.h.s. of
(\ref{msconstructions}), 
\begin{equation}
\mathfrak{m}_{3}\left(\mathbb{A}_{1},\mathbb{A}_{2}\right)=\mathfrak{der}(\mathbb{A}_{1})\oplus\mathfrak{der}(\mathbb{A}_{2})\,\dot{+}\,\mathfrak{sa}(3,\mathbb{A}_{1}\otimes\mathbb{A}_{2}).\label{Vin}
\end{equation}
This formula is manifestly symmetric under the exchange $\mathbb{A}_{1}\leftrightarrow\mathbb{A}_{2}$,
but it only holds if $\mathbb{A}_{1}$ and $\mathbb{A}_{2}$ are both
composition algebras; it involves only the derivation Lie algebras
of $\mathbb{A}$ and $\mathbb{B}$, as well as the $3\times3$ antisymmetric
matrices over $\mathbb{A}_{1}\otimes\mathbb{A}_{2}$. As mentioned
above, a suitable deformation of (\ref{Vin}) is possible, in order
to give rise to $\mathfrak{m}_{1,2}\left(\mathbb{A}_{1},\mathbb{A}_{2}\right)$.
\item Another symmetric formula under the exchange $\mathbb{A}_{1}\leftrightarrow\mathbb{A}_{2}$
was obtained by Barton and Sudbery \cite{BS} (also cf. \cite{Evans:2009ed}),
and it is given by the third line of the r.h.s. of (\ref{msconstructions}),
\begin{equation}
\mathfrak{m}_{3}\left(\mathbb{A}_{1},\mathbb{A}_{2}\right)=\mathfrak{tri}(\mathbb{A}_{1})\oplus\mathfrak{tri}(\mathbb{A}_{2})\,\dot{+}\,3(\mathbb{A}_{1}\otimes\mathbb{A}_{2}),\label{BaSu}
\end{equation}
involving the triality algebras $\mathfrak{tri}\left(\mathbb{A}_{1}\right)$
and $\mathfrak{tri}\left(\mathbb{A}_{2}\right)$, together with three
copies of the tensor product $\mathbb{A}_{1}\otimes\mathbb{A}_{2}$.
Again, as mentioned above, a suitable deformation of (\ref{BaSu})
is possible, in order to give rise to $\mathfrak{m}_{1,2}\left(\mathbb{A}_{1},\mathbb{A}_{2}\right)$. 
\end{enumerate}
Other, different versions of the construction of Freudenthal-Tits
Magic Square have been developed by Santander and Herranz \cite{Santander,Santander-2},
Atsuyama \cite{Atsuyama1,Atsuyama} and Elduque \cite{EldMS1,EldMS2,ElDuque-Comp},
all involving composition algebras (even though Elduque's construction
involves \textit{flexible} composition algebras instead of \textit{alternative}
composition algebras \cite{ElDuque-Comp}).

\subsection{\label{Order 2}Magic Square of order 2}

Inspired by the works of Freudenthal and Tits, Barton and Sudbery
\cite{BS} also considered a different Magic Square, based on $2\times2$
matrices, and defined by the following formula : 
\begin{equation}
\mathfrak{m}_{2}\left(\mathbb{A}_{1},\mathbb{A}_{2}\right):=\mathfrak{so}\left(\mathbb{A}_{1}^{\prime}\right)\oplus\mathfrak{der}\left(\mathfrak{J}_{2}\left(\mathbb{A}_{2}\right)\right)\oplus\left(\mathbb{A}_{1}^{\prime}\otimes\mathfrak{J}_{2}^{\prime}\left(\mathbb{A}_{2}\right)\right),\label{order2}
\end{equation}
which also enjoys a ``Vinberg-like'' equivalent version as 
\begin{equation}
\mathfrak{m}_{2}\left(\mathbb{A}_{1},\mathbb{A}_{2}\right)=\mathfrak{so}\left(\mathbb{A}_{1}^{\prime}\right)\oplus\mathfrak{so}\left(\mathbb{A}_{2}^{\prime}\right)\oplus\mathfrak{sa}_{2}\left(\mathbb{A}_{1}\otimes\mathbb{A}_{2}\right),\label{order2-Vin}
\end{equation}
that was more recentely used at Lie group level in \cite{Dray-Huerta-Kincaid}.
Formula (\ref{order2}) can be generalized to involve simple, rank-2
\textit{Lorentzian} Jordan algebras $\mathfrak{J}_{1,1}\left(\mathbb{A}_{2}\right)$
(see Sec. \ref{Eucl-Lor}), thus obtaining the Lorentzian version
of the Magic Square of order 2,
\begin{equation}
\mathfrak{m}_{1,1}\left(\mathbb{A}_{1},\mathbb{A}_{2}\right):=\mathfrak{so}\left(\mathbb{A}_{1}^{\prime}\right)\oplus\mathfrak{der}\left(\mathfrak{J}_{1,1}\left(\mathbb{A}_{2}\right)\right)\oplus\left(\mathbb{A}_{1}^{\prime}\otimes\mathfrak{J}_{1,1}^{\prime}\left(\mathbb{A}_{2}\right)\right),\label{order2-Lor}
\end{equation}
which we will explicitly evaluate in the treatment below (for the
first time in literature, to the best of our knowledge).

\subsection{\label{Expl}Octonionic entries of Magic Squares of order 2 and 3}

In order to rigorously define all possible octonionic \textit{``Rosenfeld
planes''} and \textit{``Rosenfeld lines''} over $\mathbb{R}$ (i.e.
all possible real forms thereof), we need to isolate all octonionic
and split-octonionic entries of all real forms of the Freudenthal-Tits
Magic Square (order 3) \cite{Squaring-Magic,Kind-of-Magic,Magic-Pyramids,BS}
and of the Magic Square of order 2 \cite{BS}, i.e. we need to consider
the entries $\mathfrak{m}_{\alpha}\left(\mathbb{A}_{1},\mathbb{A}_{2}\right)$
in which \textit{at least one} of $\mathbb{A}_{1}$ and $\mathbb{A}_{2}$
is $\mathbb{O}$ or $\mathbb{O}_{s}$, for all possible real forms,
namely for $\alpha=3$ (Euclidean order 3), $1,2$ (Lorentzian order
3), $2$ (Euclidean order 2) and $1,1$ (Lorentzian order 2).

The octonionic entries of the Euclidean Magic Squares of order 3 $\mathfrak{m}_{3}$
(i.e., of all Euclidean real forms of Freudenthal-Tits Magic Square)
are

\medskip{}

\begin{center}
\begin{tabular}{|c|c|c|c|c|}
\hline 
$\mathbb{A}$ & $\mathfrak{m}_{3}\left(\mathbb{A},\mathbb{O}\right)$ & $\mathfrak{m}_{3}\left(\mathbb{A}_{s},\mathbb{O}\right)$ & $\mathfrak{m}_{3}\left(\mathbb{O}_{s},\mathbb{A}\right)$ & $\mathfrak{m}_{3}\left(\mathbb{A}_{s},\mathbb{O}_{s}\right)$\tabularnewline
\hline 
\hline 
$\mathbb{R}$ & $\mathfrak{f}_{4\left(-52\right)}$ & $\mathfrak{f}_{4\left(-52\right)}$ & $\mathfrak{f}_{4\left(4\right)}$ & $\mathfrak{f}_{4\left(4\right)}$\tabularnewline
\hline 
$\mathbb{C}$ & $\mathfrak{e}_{6\left(-78\right)}$ & $\mathfrak{e}_{6\left(-26\right)}$ & $\mathfrak{e}_{6\left(2\right)}$ & $\mathfrak{e}_{6\left(6\right)}$\tabularnewline
\hline 
$\mathbb{H}$ & $\mathfrak{e}_{7\left(-133\right)}$ & $\mathfrak{e}_{7\left(-25\right)}$ & $\mathfrak{e}_{7\left(-5\right)}$ & $\mathfrak{e}_{7\left(7\right)}$\tabularnewline
\hline 
$\mathbb{O}$ & $\mathfrak{e}_{8\left(-248\right)}$ & $\mathfrak{e}_{8\left(-24\right)}$ & $\mathfrak{e}_{8\left(-24\right)}$ & $\mathfrak{e}_{8\left(8\right)}$\tabularnewline
\hline 
\end{tabular} 
\par\end{center}

\medskip{}

The entries of the table above comprise all real forms of $\mathfrak{e}_{8}$
and $\mathfrak{e}_{7}$, and most of the real forms of $\mathfrak{e}_{6}$
and $\mathfrak{f}_{4}$. The missing real forms $\mathfrak{f}_{4\left(-20\right)}$
and $\mathfrak{e}_{6\left(-14\right)}$ can be recovered as the octonionic
entries of the Lorentzian Magic Squares of order 3 $\mathfrak{m}_{2,1}$
(i.e., of all Lorentzian real forms of Freudenthal-Tits Magic Square),
given by

\medskip{}

\begin{center}
\begin{tabular}{|c|c|c|c|c|}
\hline 
$\mathbb{A}$ & $\mathcal{\mathfrak{m}}_{1,2}\left(\mathbb{A},\mathbb{O}\right)$ & $\mathcal{\mathfrak{m}}_{1,2}\left(\mathbb{A}_{s},\mathbb{O}\right)$ & $\mathcal{\mathfrak{m}}_{1,2}\left(\mathbb{O}_{s},\mathbb{A}\right)$ & $\mathcal{\mathfrak{m}}_{1,2}\left(\mathbb{A}_{s},\mathbb{O}_{s}\right)$\tabularnewline
\hline 
\hline 
$\mathbb{R}$ & $\mathfrak{f}_{4\left(-20\right)}$ & $\mathfrak{f}_{4\left(-20\right)}$ & $\mathfrak{f}_{4\left(4\right)}$ & $\mathfrak{f}_{4\left(4\right)}$\tabularnewline
\hline 
$\mathbb{C}$ & $\mathfrak{e}_{6\left(-14\right)}$ & $\mathfrak{e}_{6\left(-26\right)}$ & $\mathfrak{e}_{6\left(2\right)}$ & $\mathfrak{e}_{6\left(6\right)}$\tabularnewline
\hline 
$\mathbb{H}$ & $\mathfrak{e}_{7\left(-5\right)}$ & $\mathfrak{e}_{7\left(-25\right)}$ & $\mathfrak{e}_{7\left(-5\right)}$ & $\mathfrak{e}_{7\left(7\right)}$\tabularnewline
\hline 
$\mathbb{O}$ & $\mathfrak{e}_{8\left(8\right)}$ & $\mathfrak{e}_{8\left(-24\right)}$ & $\mathfrak{e}_{8\left(-24\right)}$ & $\mathfrak{e}_{8\left(8\right)}$\tabularnewline
\hline 
\end{tabular}
\par\end{center}

\medskip{}

On the other hand, the octonionic entries of the Euclidean Magic Squares
of order 2 $\mathfrak{m}_{2}$ are\footnote{The explicit form of $\mathfrak{m}_{2}\left(\mathbb{A}_{s},\mathbb{B}_{s}\right)$
is, as far as we know, not present in the current literature. We will
present a detailed treatment elsewhere, and here we confine ourselves
to report its octonionic column only.}

\medskip{}

\begin{center}
\begin{tabular}{|c|c|c|c|c|}
\hline 
$\mathbb{A}$ & $\mathcal{\mathfrak{m}}_{2}\left(\mathbb{A},\mathbb{O}\right)$ & $\mathcal{\mathfrak{m}}_{2}\left(\mathbb{A}_{s},\mathbb{O}\right)$ & $\mathcal{\mathfrak{m}}_{2}\left(\mathbb{O}_{s},\mathbb{A}\right)$ & $\mathcal{\mathfrak{m}}_{2}\left(\mathbb{A}_{s},\mathbb{O}_{s}\right)$\tabularnewline
\hline 
\hline 
$\mathbb{R}$ & $\mathfrak{so}\left(9\right)$ & $\mathfrak{so}\left(9\right)$ & $\mathfrak{so}\left(5,4\right)$ & $\mathfrak{so}\left(5,4\right)$\tabularnewline
\hline 
$\mathbb{C}$ & $\mathfrak{so}\left(10\right)$ & $\mathfrak{so}\left(9,1\right)$ & $\mathfrak{so}\left(6,4\right)$ & $\mathfrak{so}\left(5,5\right)$\tabularnewline
\hline 
$\mathbb{H}$ & $\mathfrak{so}\left(12\right)$ & $\mathfrak{so}\left(10,2\right)$ & $\mathfrak{so}\left(8,4\right)$ & $\mathfrak{so}\left(6,6\right)$\tabularnewline
\hline 
$\mathbb{O}$ & $\mathfrak{so}\left(16\right)$ & $\mathfrak{so}\left(12,4\right)$ & $\mathfrak{so}\left(12,4\right)$ & $\mathfrak{so}\left(8,8\right)$\tabularnewline
\hline 
\end{tabular}
\par\end{center}

\medskip{}

Finally, the octonionic entries of the Lorentzian Magic Squares of
order 2 $\mathfrak{m}_{1,1}$ are given by\footnote{The Lorentzian magic square of order two yield to three algebras,
i.e. $\mathfrak{so}\left(8,1\right)$, $\mathfrak{so}\left(8,2\right)$
and $\mathfrak{so}\left(8,4\right)$ which are not covered in other
magic squares. The explicit forms of $\mathfrak{m}_{1,1}\left(\mathbb{A},\mathbb{B}\right)$,
$\mathfrak{m}_{1,1}\left(\mathbb{A}_{s},\mathbb{B}\right)$ and $\mathfrak{m}_{1,1}\left(\mathbb{A}_{s},\mathbb{B}_{s}\right)$
are, as far as we know, not present in the current literature. We
will present a detailed treatment elsewhere, and here we confine ourselves
to report their octonionic rows and columns only.}

\medskip{}

\begin{center}
\begin{tabular}{|c|c|c|c|c|}
\hline 
$\mathbb{A}$ & $\mathcal{\mathcal{\mathfrak{m}}}_{1,1}\left(\mathbb{A},\mathbb{O}\right)$ & $\mathcal{\mathcal{\mathfrak{m}}}_{1,1}\left(\mathbb{A}_{s},\mathbb{O}\right)$ & $\mathcal{\mathcal{\mathfrak{m}}}_{1,1}\left(\mathbb{O}_{s},\mathbb{A}\right)$ & $\mathcal{\mathcal{\mathfrak{m}}}_{1,1}\left(\mathbb{A}_{s},\mathbb{O}_{s}\right)$\tabularnewline
\hline 
\hline 
$\mathbb{R}$ & $\mathfrak{so}\left(8,1\right)$ & $\mathfrak{so}\left(8,1\right)$ & $\mathfrak{so}\left(5,4\right)$ & $\mathfrak{so}\left(5,4\right)$\tabularnewline
\hline 
$\mathbb{C}$ & $\mathfrak{so}\left(8,2\right)$ & $\mathfrak{so}\left(9,1\right)$ & $\mathfrak{so}\left(6,4\right)$ & $\mathfrak{so}\left(5,5\right)$\tabularnewline
\hline 
$\mathbb{H}$ & $\mathfrak{so}\left(8,4\right)$ & $\mathfrak{so}\left(10,2\right)$ & $\mathfrak{so}\left(8,4\right)$ & $\mathfrak{so}\left(6,6\right)$\tabularnewline
\hline 
$\mathbb{O}$ & $\mathfrak{so}\left(8,8\right)$ & $\mathfrak{so}\left(12,4\right)$ & $\mathfrak{so}\left(12,4\right)$ & $\mathfrak{so}\left(8,8\right)$\tabularnewline
\hline 
\end{tabular}
\par\end{center}

\section{\label{2}``Magic''\ formulæ\ for Rosenfeld planes}

In his study of the geometry of Lie groups \cite{RosenfeldGroup2-1},
Rosenfeld defined its ``projective''\ planes $\left(\mathbb{A}\otimes\mathbb{B}\right)P^{2}$
as the completion of some affine planes obtained as non-associative
modules over the tensor algebra $\mathbb{A}\otimes\mathbb{B}$. He
then argued the form of the matrices composing the linear transformations
associated with the Lie algebra of the collineations that preserved
the polarity, i.e. the isometries of the plane. This approach allowed
him to relate all real forms of exceptional Lie groups \cite{Rosenfeld-1993,Rosenf98}
with isometries of suitable ``projective''\ hyperbolic spaces.
As mentioned above, unfortunately the procedure carried out by Rosenfeld
was not rigorous, since many of theorems used in \cite{RosenfeldGroup2-1}
do not extend to the case of algebras that are not alternative nor
of composition, as it is in the case of many algebras listed in Table
\ref{tab:Properties-of-the}. 

We present here a general and rigorous way to define the spaces considered
by Rosenfeld, in terms of their isometry and isotropy groups, namely
using the theory of coset manifolds.

\subsection{Coset manifolds}

Since our ``magic'' formulæ define all Rosenfeld planes as coset
manifolds, it is worth briefly reviewing them, and their relation
to homogenous spaces. An \emph{homogeneous space} is a manifold on
which a Lie group acts transitively, i.e. a manifold on which is defined
an action $\rho_{g}$ from $G\times M$ in $M$ such that $\rho_{e}\left(m\right)=m$
for $e$ the identity in $G$ and $m\in M$ and for which, given any
$m,n\in M$ it exists a a not necessarely unique $g\in G$ such that
$\rho_{g}\left(m\right)=n$. Within this framework, the \emph{isotropy
group} $\text{Isot}_{m}\left(G\right)$ is the set formed by the elements
of $G$ that fix the point $m\in M$ under the action of $G$, i.e.
\begin{equation}
\text{Isot}_{m}\left(G\right)=\left\{ g\in G:\rho_{g}\left(m\right)=m\right\} .
\end{equation}
Since, by definition of group action, we have that $\rho_{e}\left(m\right)=m$
and $\rho_{gh}\left(m\right)=\rho_{g}\left(\rho_{h}\left(m\right)\right)$,
then $K=\text{Isot}_{m}\left(G\right)$ is a closed subgroup of $G$
and, moreover the natural map from the quotient space $G/K$ in $M$
given by $gK\mapsto gm$ is a diffeomorphism (see \cite{Ar} and Refs.
therein). 
\begin{table}
\centering{}%
\begin{tabular}{|c|c|c|c|c|c|}
\hline 
Algebra & Comm. & Ass. & Alter. & Flex. & Pow. Ass.\tabularnewline
\hline 
\hline 
$\mathbb{C}\otimes\mathbb{C}$ & Yes & Yes & Yes & Yes & Yes\tabularnewline
\hline 
$\mathbb{C}\otimes\mathbb{H}$ & No & Yes & Yes & Yes & Yes\tabularnewline
\hline 
$\mathbb{H}\otimes\mathbb{H}$ & No & Yes & Yes & Yes & Yes\tabularnewline
\hline 
$\mathbb{C}\otimes\mathbb{O}$ & No & No & Yes & Yes & Yes\tabularnewline
\hline 
$\mathbb{H}\otimes\mathbb{O}$ & No & No & No & No & No\tabularnewline
\hline 
$\mathbb{O}\otimes\mathbb{O}$ & No & No & No & No & No\tabularnewline
\hline 
\end{tabular}\caption{\label{tab:Properties-of-the}Properties of the algebra $\mathbb{A}\otimes\mathbb{B}$
where $\mathbb{A},\mathbb{B}$ are Hurwitz algebras. As for the property
an algebra $\mathbb{A}$ is said to be \emph{commutative} if $xy=yx$
for every $x,y\in X$; it is defined as \emph{associative} if satisfies
$x\left(yz\right)=\left(xy\right)z$; \emph{alternative} if $x\left(yx\right)=\left(xy\right)x$;
\emph{flexible} if $x\left(yy\right)=\left(xy\right)y$ and, finally,
\emph{power-associative} if $x\left(xx\right)=\left(xx\right)x$.}
\end{table}

Given a Lie group $G$ and a closed subgroup $K<G$, then the coset
space $G/K=\left\{ gK:g\in G\right\} $ is endowed with a natural
manifold structure inherited by $G$ and is, therefore, called a \emph{coset
manifold}. Notice that the action of $G$ on the coset manifold $G/K$,
given by the translation $\tau_{g}$ defined as 
\begin{equation}
\tau_{g}\left(m\right)=gm,
\end{equation}
for every $g\in G$ and $m\in G/K$, is \emph{transitive}, i.e. for
every $m,n\in G/K$ it exists a (not necessarily unique) $g\in G$
such that $n=gm$, and thus the coset manifold $G/K$ is an homogenous
space. On the other hand, since multiple groups can act transitively
on the same manifold with different isotropy groups, then an homogeneous
space can be realised in multiple way as a coset manifold. Moreover,
a close look to the definitions shows that the isotropy group $\text{Isot}_{m}\left(G/K\right)$
is exactly the closed subgroup $K$.\footnote{In particular, the origin of $G/K$ is, by definition, the point at
which the $K$-invariance is immediately manifest.} In general, the \textit{holonomy} \textit{subgroup} and the isotropy
subgroup have the same identity-connected component; so, if one assumes
that $G/K$ is simply-connected, they are equal (see e.g. \cite{Besse,Helgason}
and Refs. therein).

Moreover, let $G$ be a Lie group and $K$ a closed and connected
subgroup of $G$, denoting with $\mathfrak{g}$ and $\mathfrak{k}$
their respective Lie algebras, then the coset manifold $G/K$ is \emph{reductive}
if there exists a subspace $\mathfrak{m}$ such that $\mathfrak{g}=\mathfrak{k}\oplus\mathfrak{m}$
and 
\begin{align}
\begin{cases}
\left[\mathfrak{k},\mathfrak{k}\right] & \subset\mathfrak{k},\\
\left[\mathfrak{k},\mathfrak{m}\right] & \subset\mathfrak{m},
\end{cases}\label{eq:homog-reduct-space}
\end{align}
while if in addition to (\ref{eq:homog-reduct-space}) we also have
\begin{equation}
\left[\mathfrak{m},\mathfrak{m}\right]\subset\mathfrak{k},\label{eq:symmetric space}
\end{equation}
then the space is \emph{symmetric}. All Rosenfeld planes are symmetric
coset manifolds. 

It is also worth noting that for any coset manifold $G/K$ the structure
constants of the Lie algebra $\mathfrak{g}$ of the Lie group $G$
define completely the structure constants of the manifold, thus the
invariant metrics, and all the metric-dependent tensors, such as the
curvature tensor, the Ricci tensor, etc. Indeed, let $\left\{ E_{1},...,E_{n}\right\} $
be a basis for $\mathfrak{g}$ in such a way that $\left\{ E_{1},...,E_{m}\right\} $
are a basis for $\mathfrak{k}$, which we will also call $\left\{ K_{1},...,K_{m}\right\} $
for readibility reasons, and $\left\{ E_{m+1},...,E_{n}\right\} $
a basis for $\mathfrak{m}$ that we will also denote as $\left\{ M_{m+1},...,M_{n}\right\} $.
Then consider the structure constants of the algebra $\mathfrak{g}$,
i.e.

\begin{equation}
\left[E_{j},E_{k}\right]=\stackrel[i=1]{n}{\sum}C_{jk}^{i}E_{i},
\end{equation}
 for every $j,k\in\left\{ 1,...,n\right\} $. Conditions for reductivity
in (\ref{eq:homog-reduct-space}) are then translated in 
\begin{equation}
\begin{cases}
\left[K_{j},K_{k}\right]=\stackrel[i=1]{m}{\sum}C_{jk}^{i}K_{i} & \text{for }j,k\in\left\{ 1,...,m\right\} ,\\
\left[K_{j},M_{k}\right]=\stackrel[i=m+1]{n}{\sum}C_{jk}^{i}M_{i} & \text{for }j\in\left\{ 1,...,m\right\} ,k\in\left\{ m+1,...,n\right\} ,
\end{cases}
\end{equation}
 while, on the other hand, from

\begin{equation}
\left[M_{j},M_{k}\right]=\stackrel[i=1]{m}{\sum}C_{jk}^{i}K_{i}+\stackrel[i=m+1]{n}{\sum}C_{jk}^{i}M_{i},
\end{equation}
we deduce that if $C_{jk}^{i}=0$ for all $i,j,k\in\left\{ m+1,...,n\right\} $,
we have also a symmetric space. From the structure constants $C_{jk}^{i}$,
all geometrical invariants of the coset manifold can be obtained (see
\cite{Fre-Fedotov} for all technical details) such as the Riemann
tensor over the coset manifold $G/K$ which is given by
\begin{equation}
R_{bcd}^{a}=\stackrel[e=m+1]{n}{\sum}\frac{1}{\lambda^{2}}\left(\frac{1}{8}C_{ed}^{a}C_{bc}^{e}-\frac{1}{8\lambda}C_{be}^{a}C_{cd}^{e}-C_{ec}^{a}C_{bd}^{e}\right)-\frac{1}{2\lambda^{2}}\stackrel[i=1]{m}{\sum}C_{bi}^{a}C_{cd}^{i},
\end{equation}
and that, in case of symmetric spaces, i.e. when $C_{jk}^{i}=0$ for
$i,j,k\in\left\{ m+1,...,n\right\} $, reduces drastically to 
\begin{equation}
R_{bcd}^{a}=-\frac{1}{2\lambda^{2}}\stackrel[i=1]{m}{\sum}C_{bi}^{a}C_{cd}^{i},
\end{equation}
for every $a,b,c,d\in\left\{ m+1,...,n\right\} $. 

In fact, what is relevant for our purposes is that the Lie group $G$,
which acts as isometry group, and its closed subgroup $K$, which
acts as isotropy group, define completely all the metric properties
and geometric invariants of the Rosenfeld planes. On the other hand,
instead of working on Lie groups $G$ and $K$, i.e. with the isometry
group and the isotropy group respectively, in the following sections
we will work with their respective Lie algebras $\mathfrak{g}$ and
$\mathfrak{k}$, in order to recover the appropriate $G$ and $K$
and thus define the Rosenfeld plane.

\subsection{``Magic''\ formulæ}

As seen in the previous section, the geometry of a coset manifold
is fully determined by the (global) isotropy Lie group $G$ and by
its (local) subgroup given by the isotropy Lie group $K$. We now
characterize the Rosenfeld ``projective''\ planes by specifying
the real form of the Lie groups that will act as isometry and isotropy
groups, while the Lie brackets - and therefore the metrical properties
of the space - are those arising from Tits construction. It is here
worth noting that the classification of the real forms of simple Lie
groups makes use of the \textit{character} $\chi$, defined as the
difference of the cardinality of non-compact and compact generators,
i.e. $\chi:=\#_{nc}-\#_{c}$; consequently, Rosenfeld planes will
have a corresponding $\chi$.\medskip{}

Let $\mathbb{A}$ and $\mathbb{B}$ be two Hurwitz algebras, in their
division or split versions. Moreover, let the Lie group $\mathcal{A}\left(\mathbb{A}\right)$
be such that its Lie algebra is Lie$\left(\mathcal{A}\left(\mathbb{A}\right)\right)\equiv\mathfrak{A}(\mathbb{A}):=\mathfrak{tri}\left(\mathbb{A}\right)\ominus\mathfrak{so}\left(\mathbb{A}\right)$,
cfr. Table \ref{tab:Triality-algebras}. Moreover, let $\mathcal{M}_{\alpha}\left(\mathbb{A},\mathbb{B}\right)$
be the Lie group with Lie algebra given by the $\left(\mathbb{A},\mathbb{B}\right)$-entry
of the Magic Square $\mathfrak{m}_{\alpha}\left(\mathbb{A},\mathbb{B}\right)$,
namely\footnote{For the pseudo-orthogonal Lie algebras, we will generally consider
the spin covering of the corresponding Lie group.} Lie$\left(\mathcal{M}_{\alpha}\left(\mathbb{A},\mathbb{B}\right)\right)=\mathfrak{m}_{\alpha}\left(\mathbb{A},\mathbb{B}\right)$.
In order to characterize ``projective''\ Rosenfeld planes over
tensor products of Hurwitz algebras in terms of coset manifolds with
isometry (resp. isotropy) Lie groups whose Lie algebras are entries
of real forms of the Magic Square of order 3 (resp. 2), we now introduce
the following three different classes of locally symmetric, (pseudo-)Riemannian
coset manifolds, that we name as \textit{Rosenfeld planes} :
\begin{enumerate}
\item The \textit{projective Rosenfeld plane} 
\begin{equation}
\left(\mathbb{A}\otimes\mathbb{B}\right)P^{2}\simeq\frac{\mathcal{M}_{3}\left(\mathbb{A},\mathbb{B}\right)}{\mathcal{M}_{2}\left(\mathbb{A},\mathbb{B}\right)\otimes\mathcal{A}\left(\mathbb{A}\right)\otimes\mathcal{A}\left(\mathbb{B}\right)}.\label{eq:projectiv}
\end{equation}
\item The \textit{hyperbolic Rosenfeld plane} 
\begin{equation}
\left(\mathbb{A}\otimes\mathbb{B}\right)H^{2}\simeq\frac{\mathcal{M}_{1,2}\left(\mathbb{A},\mathbb{B}\right)}{\mathcal{M}_{2}\left(\mathbb{A},\mathbb{B}\right)\otimes\mathcal{A}\left(\mathbb{A}\right)\otimes\mathcal{A}\left(\mathbb{B}\right)}.\label{eq:hyperb}
\end{equation}
\item The \textit{pseudo-Rosenfeld plane}
\begin{equation}
\left(\mathbb{A}\otimes\mathbb{B}\right)\widetilde{H}^{2}\simeq\frac{\mathcal{M}_{1,2}\left(\mathbb{A},\mathbb{B}\right)}{\mathcal{M}_{1,1}\left(\mathbb{A},\mathbb{B}\right)\otimes\mathcal{A}\left(\mathbb{A}\right)\otimes\mathcal{A}\left(\mathbb{B}\right)}.\label{eq:mixed}
\end{equation}
\end{enumerate}
Eqs. (\ref{eq:projectiv}), (\ref{eq:hyperb}) and (\ref{eq:mixed})
are named \textit{``magic''\ formulæ}, since they characterize
the Rosenfeld planes as homogeneous (symmetric) manifolds, with isometry
(resp. isotropy) groups whose Lie algebras are given by the entries
of some real forms of the Magic Square of order 3 (resp. 2), with
further isotropy factors given by the Lie groups $\mathcal{A}$ associated
to the algebras $\mathbb{A}$ and $\mathbb{B}$ defining the tensor
product associated to the class of Rosenfeld plane under consideration.
As previously noticed, the term ``projective'' and ``hyperbolic''
are here to be intended in a vague sense since none of the octonionic
Rosenfeld planes with dimension greater than 16 satisfy axioms of
projective or hyperbolic geometry. Nevertheless, those adjectives
are not arbitrary since the notion of such projective planes can be
made precise as in \cite{Atsuyama,Atiyah-Bernd} and it then agrees
with the above definitions. 

\section{\label{3}\textit{Octonionic} Rosenfeld planes}

We now consider the ``magic''\ formulæ\ (\ref{eq:projectiv})-(\ref{eq:mixed})
in the cases\footnote{As mentioned above, this restriction does not imply any loss of generality,
as far as the other, non-octonionic Rosenfeld planes can be obtained
as suitable sub-manifolds of the octonionic Rosenfeld planes.} in which $\mathbb{A}$ and/or $\mathbb{B}$ is $\mathbb{O}$ or $\mathbb{O}_{s}$
: this will allow us to rigorously introduce the \textit{octonionic
Rosenfeld planes}, which all share the fact that their isometry group
is a real form of an exceptional Lie group of F- or E- type; however,
it is here worth anticipating that a few real forms of Rosenfeld planes
with exceptional isometry groups cannot be characterized in this way.

\subsection{$\mathbb{O}\simeq\mathbb{R}\otimes\mathbb{O}$}

The simplest case concerns the tensor product $\mathbb{R}\otimes\mathbb{O}$,
which is nothing but $\mathbb{O}$ : in fact, this was the original
observation by Rosenfeld that started it all, yielding to the usual
octonionic projective, split-octonionic and hyperbolic plane, i.e.
$\mathbb{O}P^{2}$, $\mathbb{O}H^{2}$ and $\mathbb{O}_{s}P^{2}$,
respectively. Within the framework introduced above, the starting
point is given by the \textit{Cayley-Moufang plane} over $\mathbb{C}$,
namely by the \textit{octonionic projective plane} over $\mathbb{C}$,
i.e. by the locally symmetric coset manifold having as isometry group
the complex form of $\text{F}_{4}$, and as isotropy group $\text{Spin}\left(9,\mathbb{C}\right)$,
i.e. 
\begin{equation}
\mathbb{O}P_{\mathbb{C}}^{2}\simeq\frac{\text{F}_{4}^{\mathbb{C}}}{\text{Spin}\left(9,\mathbb{C}\right)}.\label{OP^2(C)}
\end{equation}
In this coset space formulation, the tangent space of $\mathbb{O}P_{\mathbb{C}}^{2}$
can be identified with the $\boldsymbol{16}_{\mathbb{C}}$-dimensional
spinor representation space of the isotropy group $\text{Spin}\left(9,\mathbb{C}\right)$.

Then, by specifying the formulæ\ (\ref{eq:projectiv}), (\ref{eq:hyperb})
and (\ref{eq:mixed}) for $\mathbb{A}=\mathbb{R}$ and $\mathbb{B}=\mathbb{O}$
or $\mathbb{O}_{s}$, we obtain \textit{all} real forms\footnote{By real forms of $\mathbb{O}P_{\mathbb{C}}^{2}$, we here mean the
cosets with isometry groups given by \textit{all} real (compact and
non-compact) forms of $F_{4}$, and with isotropy group given by all
(compact and non-compact) real forms of Spin$(9)$ which are subgroups
of the corresponding real form of $F_{4}$.} of (\ref{OP^2(C)}) as Rosenfeld planes over $\mathbb{R}\otimes\mathbb{O}\simeq\mathbb{O}$
or over $\mathbb{R}\otimes\mathbb{O}_{s}\simeq\mathbb{O}_{s}$; they
are summarized by the following\footnote{Recall that $\mathfrak{A}(\mathbb{R})=\mathfrak{A}(\mathbb{O})=\mathfrak{A}(\mathbb{O}_{s})=\varnothing$.}
Table \cite{Corr-Notes-Octo} :
\begin{center}
\begin{tabular}{|c|c|c|c|c|c|}
\hline 
\textbf{Plane} & \textbf{Isometry} & \textbf{Isotropy} & $\#_{nc}$ & $\#_{c}$ & $\chi$\tabularnewline
\hline 
\hline 
$\mathbb{O}P^{2}$ & $\text{F}_{4(-52)}$ & $\text{Spin}\left(9\right)$ & $0$ & $16$ & $-16$\tabularnewline
\hline 
$\mathbb{O}H^{2}$ & $\text{F}_{4(-20)}$ & $\text{Spin}\left(9\right)$ & $16$ & $0$ & $16$\tabularnewline
\hline 
$\mathbb{O}\widetilde{H}^{2}$ & $\text{F}_{4(-20)}$ & $\text{Spin}\left(8,1\right)$ & $8$ & $8$ & $0$\tabularnewline
\hline 
$\mathbb{O}_{s}P^{2}$ & $\text{F}_{4(4)}$ & $\text{Spin}\left(5,4\right)$ & $8$ & $8$ & $0$\tabularnewline
\hline 
\end{tabular}
\par\end{center}

\bigskip{}
 along with 
\begin{equation}
\mathbb{O}_{s}P^{2}\simeq\mathbb{O}_{s}H^{2}\simeq\mathbb{O}_{s}\widetilde{H}^{2}.
\end{equation}

\subsection{$\mathbb{C}\otimes\mathbb{O}$}

In the bioctonionic case, i.e. for $\mathbb{C}\otimes\mathbb{O}$,
the starting point is given by the \textit{bioctonionic projective
plane} over $\mathbb{C}$, i.e. by the locally symmetric coset manifold
having as isometry group the complex form of $E_{6}$, and as isotropy
group Spin$\left(10,\mathbb{C}\right)\otimes$U$_{1}$, i.e. 
\begin{equation}
\left(\mathbb{C}\otimes\mathbb{O}\right)P_{\mathbb{C}}^{2}\simeq\frac{\text{E}_{6}^{\mathbb{C}}}{\text{Spin}\left(10,\mathbb{C}\right)\otimes\left(\text{U}_{1}\right)_{\mathbb{C}}},\label{(CxO)P^2(C)}
\end{equation}
which is a Kähler manifold, and has been recently treated in \cite{RealF}.
In this case, the tangent space of the coset manifold (\ref{(CxO)P^2(C)})
is the $\left(\mathbf{16}_{\mathbb{C},+}\oplus\overline{\mathbf{16}}_{\mathbb{C},-}\right)$
representation of the isotropy group $\text{Spin}\left(10,\mathbb{C}\right)\otimes\text{U}_{1}$.

Then, by specifying the formulæ\ (\ref{eq:projectiv}), (\ref{eq:hyperb})
and (\ref{eq:mixed}) for $\mathbb{A}=\mathbb{C}$ or $\mathbb{C}_{s}$
and $\mathbb{B}=\mathbb{O}$ or $\mathbb{O}_{s}$, we obtain all real
forms\footnote{By real forms of $\left(\mathbb{C}\otimes\mathbb{O}\right)P_{\mathbb{C}}^{2}$,
we here mean the cosets with isometry groups given by \textit{all}
real (compact and non-compact) forms of $E_{6}$, and with isotropy
group given by all (compact and non-compact) real forms of Spin$\left(10,\mathbb{C}\right)\otimes\left(\text{U}_{1}\right)_{\mathbb{C}}$
which are subgroups of the corresponding real form of $E_{6}$.} of (\ref{(CxO)P^2(C)}) which can be expressed as Rosenfeld planes
over $\mathbb{C}$(or$~\mathbb{C}_{s}$)$\otimes\mathbb{O}$(or$~\mathbb{O}_{s}$);
they are summarized by the following\footnote{Recall that $\mathfrak{A}(\mathbb{C})=\mathfrak{u}_{1}$, and $\mathfrak{A}(\mathbb{C}_{s})=\mathfrak{so}_{1,1}$.}
Table :

\bigskip{}

\begin{center}
\begin{tabular}{|c|c|c|c|c|c|}
\hline 
\textbf{Plane} & \textbf{Isometry} & \textbf{Isotropy} & $\#_{nc}$ & $\#_{c}$ & $\chi$\tabularnewline
\hline 
\hline 
$\mathbb{\left(\mathbb{C\otimes O}\right)}P^{2}$ & $\text{E}_{6(-78)}$ & $\text{Spin}\left(10\right)\otimes\text{U}_{1}$ & $0$ & $32$ & $-32$\tabularnewline
\hline 
$\mathbb{\left(\mathbb{C\otimes O}\right)}H^{2}$ & $\text{E}_{6(-14)}$ & $\text{Spin}\left(10\right)\otimes\text{U}_{1}$ & $32$ & $0$ & $32$\tabularnewline
\hline 
$\mathbb{\left(\mathbb{C\otimes O}\right)}\widetilde{H}^{2}$ & $\text{E}_{6(-14)}$ & $\text{Spin}\left(8,2\right)\otimes\text{U}_{1}$ & $16$ & $16$ & $0$\tabularnewline
\hline 
$\left(\mathbb{C\otimes O}_{s}\right)P^{2}$ & $\text{E}_{6(2)}$ & $\text{Spin}\left(6,4\right)\otimes\text{U}_{1}$ & $16$ & $16$ & $0$\tabularnewline
\hline 
$\left(\mathbb{C}_{s}\mathbb{\otimes O}\right)P^{2}$ & $\text{E}_{6(-26)}$ & $\text{Spin}\left(9,1\right)\otimes\text{SO}\left(1,1\right)$ & $16$ & $16$ & $0$\tabularnewline
\hline 
$\left(\mathbb{C}_{s}\mathbb{\otimes O}_{s}\right)P^{2}$ & $\text{E}_{6(6)}$ & $\text{Spin}\left(5,5\right)\otimes\text{SO}\left(1,1\right)$ & $16$ & $16$ & $0$\tabularnewline
\hline 
\end{tabular}
\par\end{center}

along with
\begin{align}
\left(\mathbb{C\otimes O}_{s}\right)P^{2} & \simeq\left(\mathbb{C\otimes O}_{s}\right)H^{2}\simeq\left(\mathbb{C\otimes O}_{s}\right)\widetilde{H}^{2},\\
\left(\mathbb{C}_{s}\mathbb{\otimes O}\right)P^{2} & \simeq\left(\mathbb{C}_{s}\mathbb{\otimes O}\right)H^{2}\simeq\left(\mathbb{C}_{s}\mathbb{\otimes O}\right)\widetilde{H}^{2},\\
\left(\mathbb{C}_{s}\mathbb{\otimes O}_{s}\right)P^{2} & \simeq\left(\mathbb{C}_{s}\mathbb{\otimes O}_{s}\right)H^{2}\simeq\left(\mathbb{C}_{s}\mathbb{\otimes O}_{s}\right)\widetilde{H}^{2},
\end{align}
expressing the fact that projective, hyperbolic and pseudo Rosenfeld
planes involving $\mathbb{C}_{s}$ and/or $\mathbb{O}_{s}$ are all
isomorphic.

The first four manifolds of the above Table, having a U$_{1}$ factor
in the isotropy group, are Kähler manifolds, whereas the last two,
having a SO$(1,1)$ factor in the isotropy group, are pseudo-Kähler
manifolds.

Finally, and more importantly, there are two real forms of (\ref{(CxO)P^2(C)}),
namely the locally symmetric, pseudo-Riemannian coset Kähler manifolds
\begin{eqnarray}
 & X_{32,I} & :=\frac{\text{E}_{6(2)}}{\text{SO}^{\ast}(10)\otimes\text{U}_{1}},~\#_{nc}=20,\#_{c}=12\Rightarrow\chi=8;\label{32-1}\\
 & X_{32,II} & :=\frac{\text{E}_{6(-14)}}{\text{SO}^{\ast}(10)\otimes\text{U}_{1}},~\#_{nc}=12,\#_{c}=20\Rightarrow\chi=-8,\label{32-2}
\end{eqnarray}
whose isotropy Lie group has the corresponding Lie algebra which is
not an entry of any real form of the Magic Square of order 2.

In other words, since the Lie algebra $\mathfrak{so}^{\ast}(10)$
does not occur in any real form of the Magic Square of order 2 (see
Secs. \ref{Order 2} and \ref{Expl}), the symmetric Kähler manifolds
(\ref{32-1}) and (\ref{32-2}) cannot be characterized as Rosenfeld
planes over $\mathbb{C}$(or$~\mathbb{C}_{s}$)$\otimes\mathbb{O}$(or$~\mathbb{O}_{s}$).

\subsection{$\mathbb{H}\otimes\mathbb{O}$}

In the quaternoctonionic case, i.e. for $\mathbb{H}\otimes\mathbb{O}$,
the starting point is given by the \textit{quaternoctonionic ``projective''\ plane}
over $\mathbb{C}$, i.e. by the locally symmetric coset manifold having
as isometry group the complex form of $E_{7}$, and as isotropy group
Spin$\left(12,\mathbb{C}\right)\otimes\text{SL}\left(2,\mathbb{C}\right)$,
i.e. 
\begin{equation}
\left(\mathbb{H}\otimes\mathbb{O}\right)P_{\mathbb{C}}^{2}\simeq\frac{\text{E}_{7}^{\mathbb{C}}}{\text{Spin}\left(12,\mathbb{C}\right)\otimes\text{SL}\left(2,\mathbb{C}\right)},\label{(HxO)P^2(C)}
\end{equation}
which is a quaternionic Kähler manifold, and whose tangent space is
given by the $\left(\mathbf{32}^{(\prime)},\mathbf{2}\right)_{\mathbb{C}}$
representation\footnote{For the possible priming of the semispinor $\mathbf{32}$ of Spin$(12)$,
see e.g. \cite{Minchenko}.} of the isotropy group $\text{Spin}\left(12,\mathbb{C}\right)\otimes\text{SL}\left(2,\mathbb{C}\right)$.

Then, by specifying the formulæ\ (\ref{eq:projectiv}), (\ref{eq:hyperb})
and (\ref{eq:mixed}) for $\mathbb{A}=\mathbb{H}$ or $\mathbb{H}_{s}$
and $\mathbb{B}=\mathbb{O}$ or $\mathbb{O}_{s}$, we obtain all real
forms\footnote{By real forms of $\left(\mathbb{H}\otimes\mathbb{O}\right)P_{\mathbb{C}}^{2}$,
we here mean the cosets with isometry groups given by \textit{all}
real (compact and non-compact) forms of $E_{7}$, and with isotropy
group given by all (compact and non-compact) real forms of Spin$\left(12,\mathbb{C}\right)\otimes SL_{2}\left(\mathbb{C}\right)$
which are subgroups of the corresponding real form of $E_{7}$.} of (\ref{(HxO)P^2(C)}) which can be expressed as Rosenfeld planes
over $\mathbb{H}$(or$~\mathbb{H}_{s}$)$\otimes\mathbb{O}$(or$~\mathbb{O}_{s}$);
they are summarized by the following\footnote{Recall that $\mathfrak{A}(\mathbb{H})=\mathfrak{su}_{2}$, and $\mathfrak{A}(\mathbb{H}_{s})=\mathfrak{sl}_{2}(\mathbb{R})$.}
Table :
\begin{center}
\begin{tabular}{|c|c|c|c|c|c|}
\hline 
\textbf{Plane} & \textbf{Isometry} & \textbf{Isotropy} & $\#_{nc}$ & $\#_{c}$ & $\chi$\tabularnewline
\hline 
\hline 
$\mathbb{\left(\mathbb{H\otimes O}\right)}P^{2}$ & $\text{E}_{7(-133)}$ & $\text{Spin}\left(12\right)\otimes\text{SU}\left(2\right)$ & $0$ & $64$ & $-64$\tabularnewline
\hline 
$\mathbb{\left(\mathbb{H\otimes O}\right)}H^{2}$ & $\text{E}_{7(-5)}$ & $\text{Spin}\left(12\right)\otimes\text{SU}\left(2\right)$ & $64$ & $0$ & $64$\tabularnewline
\hline 
$\mathbb{\left(\mathbb{H\otimes O}\right)}\widetilde{H}^{2}$ & $\text{E}_{7(-5)}$ & $\text{Spin}\left(8,4\right)\otimes\text{SU}\left(2\right)$ & $32$ & $32$ & $0$\tabularnewline
\hline 
$\left(\mathbb{H}_{s}\mathbb{\otimes O}\right)P^{2}$ & $\text{E}_{7(-25)}$ & $\text{Spin}\left(10,2\right)\otimes\text{SL}\left(2,\mathbb{R}\right)$ & $32$ & $32$ & $0$\tabularnewline
\hline 
$\left(\mathbb{H}_{s}\mathbb{\otimes O}_{s}\right)P^{2}$ & $\text{E}_{7(7)}$ & $\text{Spin}\left(6,6\right)\otimes\text{SL}\left(2,\mathbb{R}\right)$ & $32$ & $32$ & $0$\tabularnewline
\hline 
\end{tabular}
\par\end{center}

\bigskip{}
 along with
\begin{align}
\mathbb{\left(\mathbb{H\otimes O}\right)}\widetilde{H}^{2} & \simeq\left(\mathbb{H\otimes O}_{s}\right)P^{2}\simeq\left(\mathbb{H\otimes O}_{s}\right)H^{2}\simeq\left(\mathbb{H}\mathbb{\otimes O}_{s}\right)\widetilde{H}^{2},\\
\left(\mathbb{H}_{s}\mathbb{\otimes O}\right)P^{2} & \cong\left(\mathbb{H}_{s}\mathbb{\otimes O}\right)H^{2}\cong\left(\mathbb{H}_{s}\mathbb{\otimes O}\right)\widetilde{H}^{2},\\
\left(\mathbb{H}_{s}\mathbb{\otimes O}_{s}\right)P^{2} & \cong\left(\mathbb{H}_{s}\mathbb{\otimes O}_{s}\right)H^{2}\cong\left(\mathbb{H}_{s}\mathbb{\otimes O}_{s}\right)\widetilde{H}^{2}.
\end{align}
The planes that do not involve split algebras, namely the ones having
a SU$(2)$ factor in the isotropy group, are quaternionic Kähler manifolds,
whereas all the ones involving split algebras, namely the ones having
a SL$(2,\mathbb{R})$ factor in the istropy group, are para-quaternionic
Kähler manifolds.

Again, there are three real forms of (\ref{(HxO)P^2(C)}), namely
the locally symmetric, pseudo-Riemannian (para-)quaternionic coset
manifolds
\begin{eqnarray}
 & X_{64,I} & :=\frac{\text{E}_{7(7)}}{\text{SO}^{\ast}\left(12\right)\otimes\text{SU}\left(2\right)},~\#_{nc}=40,\#_{c}=24\Rightarrow\chi=16;\label{64-1}\\
 & X_{64,II} & :=\frac{\text{E}_{7(-5)}}{\text{SO}^{\ast}\left(12\right)\otimes\text{SL}\left(2,\text{\ensuremath{\mathbb{R}}}\right)},~\#_{nc}=32,\#_{c}=32\Rightarrow\chi=0;\label{64-2}\\
 & X_{64,III} & :=\frac{\text{E}_{7(-25)}}{\text{SO}^{\ast}\left(12\right)\otimes\text{SU}\left(2\right)},~\#_{nc}=24,\#_{c}=40\Rightarrow\chi=-16,\label{64-3}
\end{eqnarray}
whose isotropy Lie group -up to the factor SU$\left(2\right)$ or
SL$\left(2,\text{\ensuremath{\mathbb{R}}}\right)$- has the corresponding
Lie algebra which is not an entry of any real form of the Magic Square
of order 2.

In other words, since the Lie algebra $\mathfrak{so}^{\ast}(12)$
does not occur in any real form of the Magic Square of order 2 (see
Secs. \ref{Order 2} and \ref{Expl}), the symmetric (para-)quaternionic
Kähler manifolds (\ref{64-1})-(\ref{64-3}) cannot be characterized
as Rosenfeld planes over $\mathbb{H}$(or$~\mathbb{H}_{s}$)$\otimes\mathbb{O}$(or$~\mathbb{O}_{s}$).
It should however be noticed that the isotropy group of (\ref{64-1})-(\ref{64-3})
admits an interpretation in terms of real forms of the Magic Square
of order 3, namely

\begin{eqnarray}
X_{64,I}\simeq &  & \frac{\mathcal{M}_{3}\left(\mathbb{H}_{s},\mathbb{O}_{s}\right)}{\mathcal{M}_{3}\left(\mathbb{H}_{s},\mathbb{H}\right)\otimes\mathcal{A}(\mathbb{H})},\\
X_{64,II}\simeq &  & \frac{\mathcal{M}_{3}\left(\mathbb{H},\mathbb{O}_{s}\right)}{\mathcal{M}_{3}\left(\mathbb{H},\mathbb{H}_{s}\right)\otimes\mathcal{A}(\mathbb{H}_{s})},\\
X_{64,III}\simeq &  & \frac{\mathcal{M}_{3}\left(\mathbb{H}_{s},\mathbb{O}\right)}{\mathcal{M}_{3}\left(\mathbb{H}_{s},\mathbb{H}\right)\otimes\mathcal{A}(\mathbb{H})},
\end{eqnarray}
but still the \textit{rationale} (if any) of such formulæ\ is missing,
and it is surely not the one underlying the ``magic''\ formulæ\ (\ref{eq:projectiv})-(\ref{eq:mixed}).

\subsection{$\mathbb{O}\otimes\mathbb{O}$}

In the octooctonionic case, i.e. for $\mathbb{O}\otimes\mathbb{O}$,
the starting point is given by the \textit{octoooctonionic ``projective''\ plane}
over $\mathbb{C}$, i.e. by the locally symmetric coset manifold having
as isometry group the complex form of $\text{E}_{8}$, and as isotropy
group Spin$\left(16,\mathbb{C}\right)$, i.e. 
\begin{equation}
\left(\mathbb{O}\otimes\mathbb{O}\right)P_{\mathbb{C}}^{2}\simeq\frac{\text{E}_{8}^{\mathbb{C}}}{\text{Spin}\left(16,\mathbb{C}\right)},\label{(OxO)P^2(C)}
\end{equation}
whose tangent space is given by the $\left(\boldsymbol{128}^{(\prime)}\right)_{\mathbb{C}}$
representation\footnote{Again, for the possible priming of the semispinor $\mathbf{128}$
of Spin$(16)$, see e.g. \cite{Minchenko}.} of $\text{Spin}\left(16,\mathbb{C}\right)$.

Then, by specifying the formulæ\ (\ref{eq:projectiv}), (\ref{eq:hyperb})
and (\ref{eq:mixed}) for $\mathbb{A}=\mathbb{O}$ or $\mathbb{O}_{s}$
and $\mathbb{B}=\mathbb{O}$ or $\mathbb{O}_{s}$, we obtain all real
forms\footnote{By real forms of $\left(\mathbb{O}\otimes\mathbb{O}\right)P_{\mathbb{C}}^{2}$,
we here mean the cosets with isometry groups given by \textit{all}
real (compact and non-compact) forms of $E_{8}$, and with isotropy
group given by all (compact and non-compact) real forms of Spin$\left(16,\mathbb{C}\right)$
which are subgroups of the corresponding real form of $E_{8}$.} of (\ref{(OxO)P^2(C)}) which can be expressed as Rosenfeld planes
over $\mathbb{O}$(or$~\mathbb{O}_{s}$)$\otimes\mathbb{O}$(or$~\mathbb{O}_{s}$);
they are summarized by the following Table :

\bigskip{}

\begin{center}
\begin{tabular}{|c|c|c|c|c|c|}
\hline 
\textbf{Plane} & \textbf{Isometry} & \textbf{Isotropy} & $\#_{nc}$ & $\#_{c}$ & $\chi$\tabularnewline
\hline 
\hline 
$\mathbb{\left(\mathbb{O\otimes O}\right)}P^{2}$ & $\text{E}_{8(-248)}$ & $\text{Spin}\left(16\right)$ & $0$ & $128$ & $-128$\tabularnewline
\hline 
$\mathbb{\left(\mathbb{O\otimes O}\right)}H^{2}$ & $\text{E}_{8(8)}$ & $\text{Spin}\left(16\right)$ & $128$ & $0$ & $128$\tabularnewline
\hline 
$\mathbb{\left(\mathbb{O\otimes O}\right)}\widetilde{H}^{2}$ & $\text{E}_{8(8)}$ & $\text{Spin}\left(8,8\right)$ & $64$ & $64$ & $0$\tabularnewline
\hline 
$\left(\mathbb{O}_{s}\mathbb{\otimes O}\right)P^{2}$ & $\text{E}_{8(-24)}$ & $\text{Spin}\left(12,4\right)$ & $64$ & $64$ & $0$\tabularnewline
\hline 
\end{tabular}
\par\end{center}

\bigskip{}
 along with

\begin{align}
\mathbb{\left(\mathbb{O\otimes O}\right)}\widetilde{H}^{2} & \simeq\left(\mathbb{O}_{s}\mathbb{\otimes O}_{s}\right)P^{2}\simeq\left(\mathbb{O}_{s}\mathbb{\otimes O}_{s}\right)H^{2}\simeq\left(\mathbb{O}_{s}\mathbb{\otimes O}_{s}\right)\widetilde{H}^{2},\\
\left(\mathbb{O\otimes O}_{s}\right)P^{2} & \simeq\left(\mathbb{O\otimes O}_{s}\right)H^{2}\simeq\left(\mathbb{O}\mathbb{\otimes O}_{s}\right)\widetilde{H}^{2}.
\end{align}

Again, there are two real forms of (\ref{(OxO)P^2(C)}), namely the
locally symmetric, pseudo-Riemannian coset manifolds 
\begin{equation}
X_{128,I}:=\frac{\text{E}_{8(8)}}{\text{SO}^{\ast}\left(16\right)},~\#_{nc}=72,\#_{c}=56\Rightarrow\chi=16;\label{128-1}
\end{equation}
\begin{equation}
X_{128,II}:=\frac{\text{E}_{8(-24)}}{\text{SO}^{\ast}\left(16\right)},~\#_{nc}=56,\#_{c}=72\Rightarrow\chi=-16,\label{128-2}
\end{equation}
whose isotropy Lie group has the corresponding Lie algebra which is
not an entry of any real form of the Magic Square of order 2.

In other words, since the Lie algebra $\mathfrak{so}^{\ast}(16)$
does not occur in any real form of the Magic Square of order 2 (see
Secs. \ref{Order 2} and \ref{Expl}), the symmetric manifolds (\ref{128-1})-(\ref{128-2})
cannot seemingly be characterized as Rosenfeld planes over $\mathbb{O}$(or$~\mathbb{O}_{s}$)$\otimes\mathbb{O}$(or$~\mathbb{O}_{s}$).
It should however be noticed that the isotropy group SO$^{\ast}(16)$
of (\ref{128-1})-(\ref{128-2}) admits an interpretation in terms
of the Magic Square of order 4, namely
\begin{equation}
X_{128,I}\simeq\frac{\mathcal{M}_{3}\left(\mathbb{O}_{s},\mathbb{O}_{s}\right)}{\mathcal{M}_{4}\left(\mathbb{H},\mathbb{H}_{s}\right)};
\end{equation}
\begin{equation}
X_{128,II}\simeq\frac{\mathcal{M}_{3}\left(\mathbb{O},\mathbb{O}_{s}\right)}{\mathcal{M}_{4}\left(\mathbb{H},\mathbb{H}_{s}\right)},
\end{equation}
but still the \textit{rationale} (if any) of such formulæ\ is missing,
and it is surely not the one underlying the ``magic''\ formulæ\ (\ref{eq:projectiv})-(\ref{eq:mixed}).

\begin{figure}
\includegraphics[scale=0.55]{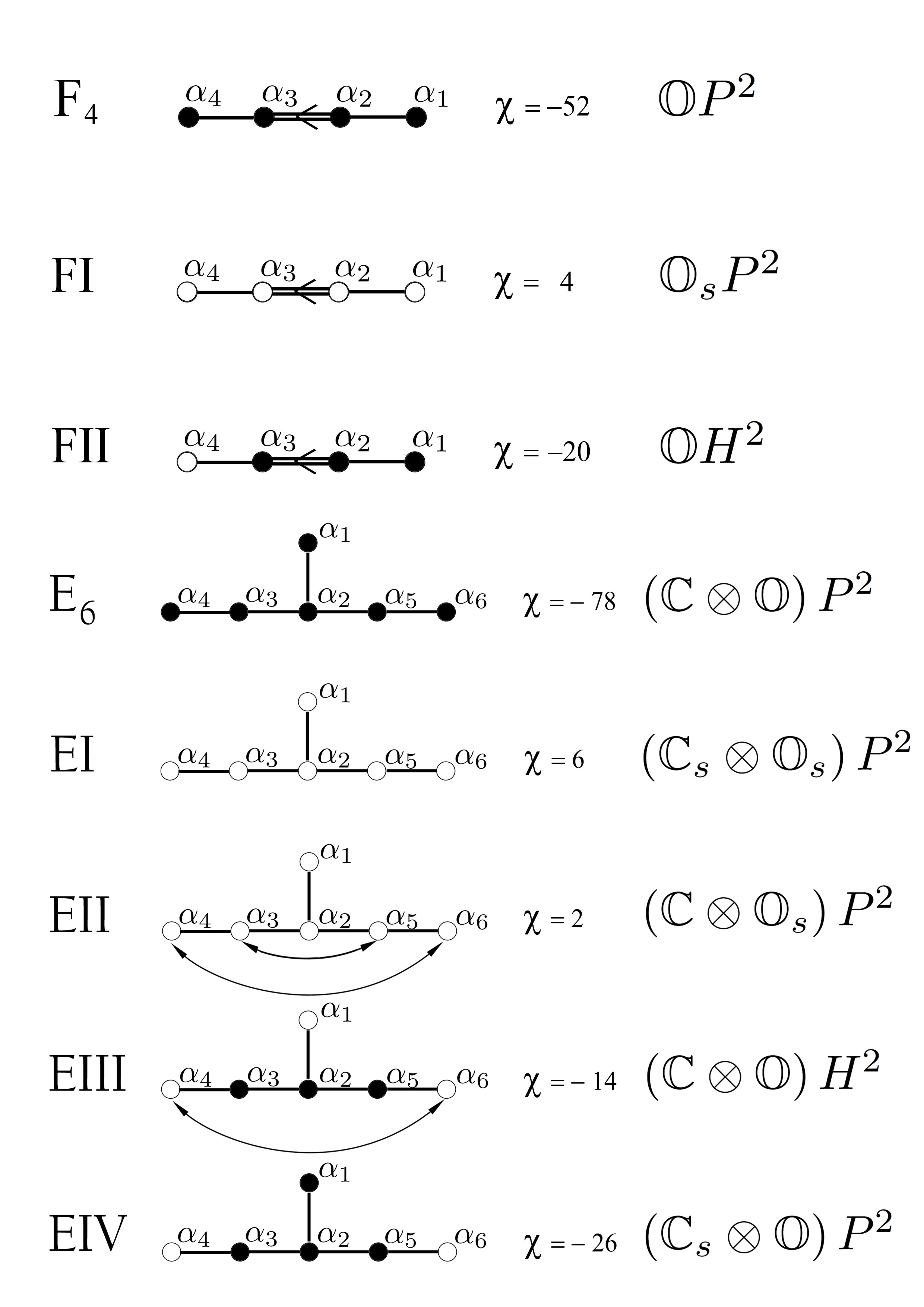}\caption{\emph{Cartan classification, Satake diagram, character $\chi$ of
exceptional Lie groups $\text{F}_{4}$ and $\text{E}_{6}$ and related
octonionic projective or hyperbolic Rosenfeld plane of which they
are isometry group.}}
\end{figure}

\begin{figure}
\includegraphics[scale=0.5]{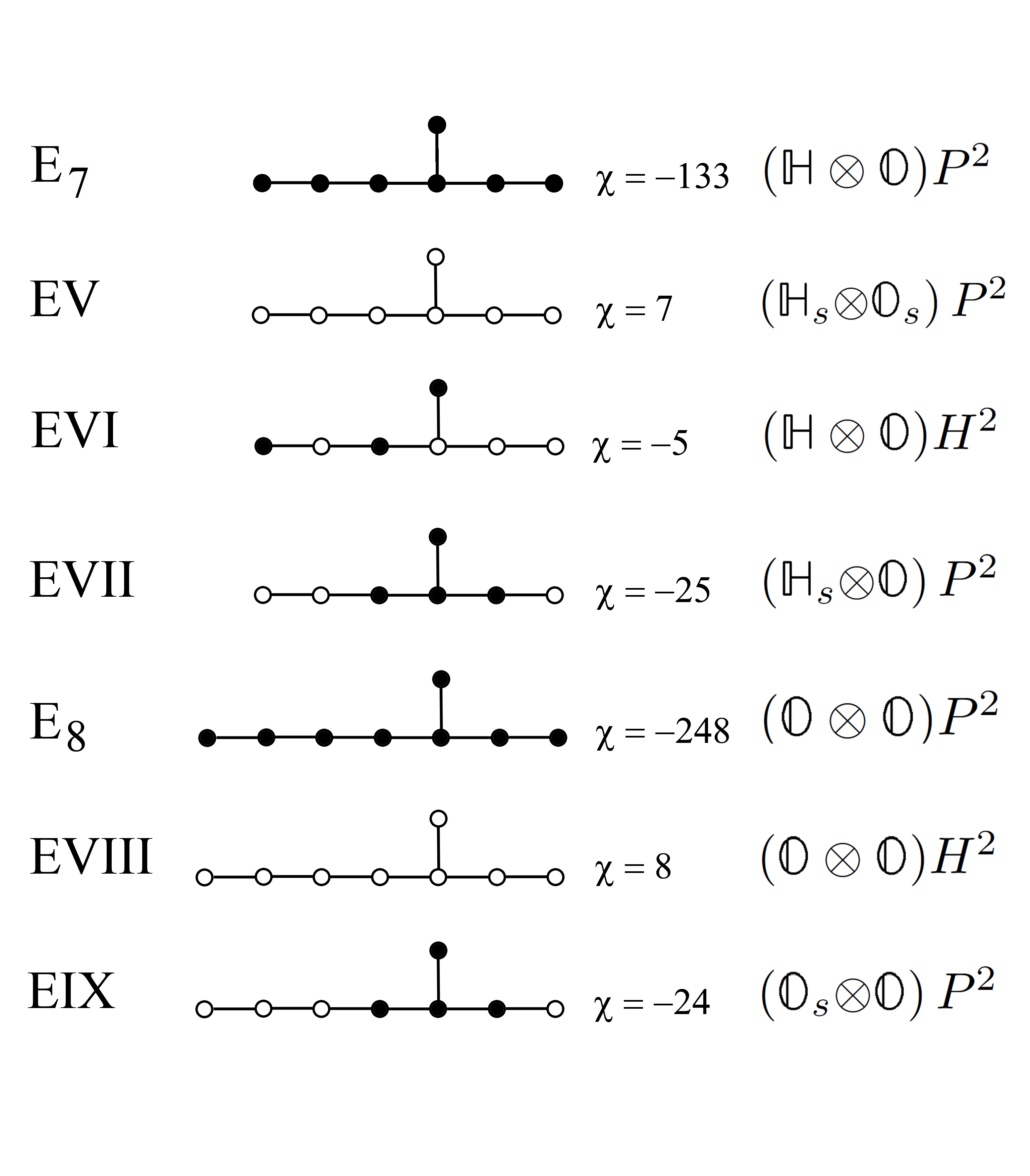}\caption{\emph{Cartan classification, Satake diagram, character $\chi$ of
exceptional Lie groups $\text{E}_{7}$ and $\text{E}_{8}$ and related
octonionic projective or hyperbolic Rosenfeld plane of which they
are isometry group.}}
\end{figure}

\section{\label{pre-4}``Magic''\ formulæ\ for Rosenfeld lines}

Let Spin$\left(\mathbb{A}\right)$ the spin covering Lie group whose
Lie algebra is Lie$\left(\text{Spin}\left(\mathbb{A}\right)\right)=\mathfrak{so}\left(\mathbb{A}\right)$,
namely the Lie algebra which preserves the norm of the Hurwitz algebra
$\mathbb{A}$. In order to characterize ``projective''\ Rosenfeld
lines over tensor products of Hurwitz algebras in terms of coset manifolds
with isometry and isotropy Lie groups whose Lie algebras respectively
are entries of real forms of the Magic Square of order 2 and norm-preserving
Lie algebras of the involved Hurwitz algebras, we now introduce a
variation of (\ref{eq:projectiv}) and (\ref{eq:hyperb}), giving
rise to the following two different classes of locally symmetric,
(pseudo-)Riemannian coset manifolds, that we name as \textit{Rosenfeld
lines} :
\begin{enumerate}
\item The ``\textit{projective''\ Rosenfeld line} 
\begin{equation}
\left(\mathbb{A}\otimes\mathbb{B}\right)P^{1}\simeq\frac{\mathcal{M}_{2}\left(\mathbb{A},\mathbb{B}\right)}{\text{Spin}\left(\mathbb{A}\right)\otimes\text{Spin}\left(\mathbb{B}\right)}.\label{eq:projectiv-line}
\end{equation}
\item The \textit{hyperbolic Rosenfeld line} 
\begin{equation}
\left(\mathbb{A}\otimes\mathbb{B}\right)H^{2}\simeq\frac{\mathcal{M}_{1,1}\left(\mathbb{A},\mathbb{B}\right)}{\text{Spin}\left(\mathbb{A}\right)\otimes\text{Spin}\left(\mathbb{B}\right)}.\label{eq:hyperb-line}
\end{equation}
\end{enumerate}
Eqs. (\ref{eq:projectiv-line}) and (\ref{eq:hyperb-line}) are named
\textit{``magic''\ formulæ}, since they characterize the \textit{Rosenfeld
lines} as homogeneous (symmetric) manifolds, with isometry and isotropy
Lie groups respectively given by the entries of some real forms of
the Magic Square of order 2 and by norm-preserving Lie algebras of
the involved Hurwitz algebras.

\section{\label{4}\textit{Octonionic} Rosenfeld lines}

We now consider the ``magic''\ formulæ\ (\ref{eq:projectiv-line})-(\ref{eq:hyperb-line})
in the cases\footnote{Also in this case, this restriction does not imply any loss of generality,
as far as the other, non-octonionic Rosenfeld lines can be obtained
as suitable sub-manifolds of the octonionic Rosenfeld lines.} in which $\mathbb{A}$ and/or $\mathbb{B}$ is $\mathbb{O}$ or $\mathbb{O}_{s}$
: this will allow us to rigorously introduce the \textit{octonionic
Rosenfeld lines}; however, again, it is here worth anticipating that
a few real forms of octonionic Rosenfeld lines cannot be characterized
in this way.

\subsection{$\mathbb{O}\simeq\mathbb{R}\otimes\mathbb{O}$}

Within the framework introduced above, the starting point is given
by the \textit{octonionic projective line} over $\mathbb{C}$, which
is nothing but the 8-sphere $S_{\mathbb{C}}^{8}$ over $\mathbb{C}$,
i.e. the locally symmetric coset manifold having as isometry group
Spin$\left(9,\mathbb{C}\right)$ and as isotropy group $\text{Spin}\left(9,\mathbb{C}\right)$
: 
\begin{equation}
\mathbb{O}P_{\mathbb{C}}^{1}\simeq\frac{\text{Spin}\left(9,\mathbb{C}\right)}{\text{Spin}\left(8,\mathbb{C}\right)}.\label{OP^1(C)}
\end{equation}
In this coset space formulation, the tangent space of $\mathbb{O}P_{\mathbb{C}}^{1}$
can be identified with the $\mathbf{8}_{(v,s,c),\mathbb{C}}$ representation
of $\text{Spin}\left(8,\mathbb{C}\right)$ where $\mathfrak{d}_{4}$-triality
\cite{D4-triality} allows to equivalently choose $v$, $s$, or $c$.

Then, by specifying the formulæ\ (\ref{eq:projectiv-line}) and (\ref{eq:hyperb-line})
for $\mathbb{A}=\mathbb{R}$ and $\mathbb{B}=\mathbb{O}$ or $\mathbb{O}_{s}$,
we obtain all real forms\footnote{By real forms of $\mathbb{O}P_{\mathbb{C}}^{1}$, we here mean the
cosets with isometry groups given by \textit{all} real (compact and
non-compact) forms of Spin$\left(9,\mathbb{C}\right)$, and with isotropy
group given by all (compact and non-compact) real forms of Spin$\left(8,\mathbb{C}\right)$
which are subgroups of the corresponding real form of Spin$\left(9,\mathbb{C}\right)$.} of (\ref{OP^1(C)}) which can be expressed as Rosenfeld lines over
$\mathbb{R}\otimes\mathbb{O}\simeq\mathbb{O}$ or over $\mathbb{R}\otimes\mathbb{O}_{s}\simeq\mathbb{O}_{s}$;
they are summarized by the following Table :

\bigskip{}

\begin{center}
\begin{tabular}{|c|c|c|c|c|c|}
\hline 
\textbf{Plane} & \textbf{Isometry} & \textbf{Isotropy} & $\#_{nc}$ & $\#_{c}$ & $\chi$\tabularnewline
\hline 
\hline 
$\mathbb{O}P^{1}$ & $\text{Spin}\left(9\right)$ & $\text{Spin}\left(8\right)$ & $0$ & $8$ & $-8$\tabularnewline
\hline 
$\mathbb{O}H^{1}$ & $\text{Spin}\left(8,1\right)$ & $\text{Spin}\left(8\right)$ & $8$ & $0$ & $8$\tabularnewline
\hline 
$\mathbb{O}_{s}P^{1}$ & $\text{Spin}\left(5,4\right)$ & $\text{Spin}\left(4,4\right)$ & $4$ & $4$ & $0$\tabularnewline
\hline 
\end{tabular} 
\par\end{center}

along with $\mathbb{O}_{s}H^{1}\simeq\mathbb{O}_{s}P^{1}$.

It should then be remarked that the octonionic line $\mathbb{O}P^{1}$
can be identified with the 8-sphere $S^{8}\equiv S_{\mathbb{R}}^{8}$
over $\mathbb{R}$, the hyperbolic line $\mathbb{O}H^{1}$ with the
8-hyperboloid $H^{8}$ and the split-octonionic line $\mathbb{O}_{s}P^{1}$
with the Kleinian 8-hyperboloid $H^{4,4}$ (for some applications
to physics and to non-compact versions of Hopf maps, see e.g. \cite{Hasebe}).

The real forms of $\mathbb{O}P_{\mathbb{C}}^{1}$ (\ref{OP^1(C)})
have the general structures
\begin{equation}
\underset{\text{structure~}I}{\frac{\text{Spin}\left(p,9-p\right)}{\text{Spin}\left(p,8-p\right)}}~\,\,\text{ or \,\,}~\underset{\text{structure~}II}{\frac{\text{Spin}\left(p,9-p\right)}{\text{Spin}\left(p-1,9-p\right)}},\label{real-O}
\end{equation}
with $p=0,1,...,9$, which gets exchanged under $p\leftrightarrow9-p$.

$\mathbb{O}P^{1}$, $\mathbb{O}H^{1}$ and $\mathbb{O}_{s}P^{1}\simeq\mathbb{O}_{s}H^{1}$
respectively correspond to structure $I$ with $p=0$, $8$, $4$
(or structure $II$ with $p=9$, $1$, $5$).

All other real forms of $\mathbb{O}P_{\mathbb{C}}^{1}$ (\ref{OP^1(C)})
cannot be characterized as Rosenfeld lines over $\mathbb{O}$ (or$~\mathbb{O}_{s}$).

\subsection{$\mathbb{C}\otimes\mathbb{O}$}

In the bioctonionic case, i.e. for $\mathbb{C}\otimes\mathbb{O}$,
the starting point is given by the \textit{bioctonionic projective
line} over $\mathbb{C}$, i.e. by the locally symmetric coset manifold
having as isometry group $\text{Spin}\left(10,\mathbb{C}\right)$
with $\text{Spin}\left(8,\mathbb{C}\right)\otimes$U$\left(1\right)$
as isotropy group : 
\begin{equation}
\left(\mathbb{C}\otimes\mathbb{O}\right)P_{\mathbb{C}}^{1}\simeq\frac{\text{Spin}\left(10,\mathbb{C}\right)}{\text{Spin}\left(8,\mathbb{C}\right)\otimes\left(\text{U}\left(1\right)\right)_{\mathbb{C}}}.\label{(CxO)P^1(C)}
\end{equation}
In this coset space formulation, the tangent space of $\left(\mathbb{C}\otimes\mathbb{O}\right)P_{\mathbb{C}}^{1}$
can be identified with the $\mathbf{8}_{(v,s,c),\mathbb{C}+}\oplus\mathbf{8}_{(v,s,c),\mathbb{C}-}$
representation of $\text{Spin}\left(8,\mathbb{C}\right)\otimes\left(\text{U}\left(1\right)\right)_{\mathbb{C}}$,
where, as above, $\mathfrak{d}_{4}$-triality \cite{D4-triality}
allows to equivalently choose $v$, $s$, or $c$.

Then, by specifying the formulæ\ (\ref{eq:projectiv-line}) and (\ref{eq:hyperb-line})
for $\mathbb{A}=\mathbb{C}$ or $\mathbb{C}_{s}$ and $\mathbb{B}=\mathbb{O}$
or $\mathbb{O}_{s}$, we obtain all real forms\footnote{By real forms of $\left(\mathbb{C}\otimes\mathbb{O}\right)P_{\mathbb{C}}^{1}$,
we here mean the cosets with isometry groups given by \textit{all}
real (compact and non-compact) forms of Spin$\left(10,\mathbb{C}\right)$,
and with isotropy group given by all (compact and non-compact) real
forms of Spin$\left(8,\mathbb{C}\right)\otimes\left(\text{U}\left(1\right)\right)_{\mathbb{C}}$
which are subgroups of the corresponding real form of Spin$\left(10,\mathbb{C}\right)$.} of (\ref{(CxO)P^1(C)}) which can be expressed as Rosenfeld lines
over $\mathbb{C}$(or$~\mathbb{C}_{s}$)$\otimes\mathbb{O}$(or$~\mathbb{O}_{s}$);
they are summarized by the following Table :

\bigskip{}

\begin{center}
\begin{tabular}{|c|c|c|c|c|c|}
\hline 
\textbf{Plane} & \textbf{Isometry} & \textbf{Isotropy} & $\#_{nc}$ & $\#_{c}$ & $\chi$\tabularnewline
\hline 
\hline 
$\left(\mathbb{C}\otimes\mathbb{O}\right)P^{1}$ & $\text{Spin}\left(10\right)$ & $\text{Spin}\left(8\right)\otimes\text{U}\left(1\right)$ & $0$ & $16$ & $-16$\tabularnewline
\hline 
$\left(\mathbb{C}\otimes\mathbb{O}\right)H^{1}$ & $\text{Spin}\left(8,2\right)$ & $\text{Spin}\left(8\right)\otimes\text{U}\left(1\right)$ & $16$ & $0$ & $16$\tabularnewline
\hline 
$\left(\mathbb{C}\otimes\mathbb{O}_{s}\right)P^{1}$ & $\text{Spin}\left(6,4\right)$ & $\text{Spin}\left(4,4\right)\otimes\text{U}\left(1\right)$ & $8$ & $8$ & $0$\tabularnewline
\hline 
$\left(\mathbb{C}_{s}\otimes\mathbb{O}\right)P^{1}$ & $\text{Spin}\left(9,1\right)$ & $\text{Spin}\left(8\right)\otimes\text{SO}\left(1,1\right)$ & $8$ & $8$ & $0$\tabularnewline
\hline 
$\left(\mathbb{C}_{s}\otimes\mathbb{O}_{s}\right)P^{1}$ & $\text{Spin}\left(5,5\right)$ & $\text{Spin}\left(4,4\right)\otimes\text{SO}\left(1,1\right)$ & $8$ & $8$ & $0$\tabularnewline
\hline 
\end{tabular} 
\par\end{center}

\bigskip{}
 along with 
\begin{align}
\left(\mathbb{C\otimes O}_{s}\right)P^{1} & \simeq\left(\mathbb{C\otimes O}_{s}\right)H^{1};\\
\left(\mathbb{C}_{s}\mathbb{\otimes O}\right)P^{1} & \simeq\left(\mathbb{C}_{s}\mathbb{\otimes O}\right)H^{1};\\
\left(\mathbb{C}_{s}\mathbb{\otimes O}_{s}\right)P^{1} & \simeq\left(\mathbb{C}_{s}\mathbb{\otimes O}_{s}\right)H^{1}.
\end{align}

The real forms of $\left(\mathbb{C}\otimes\mathbb{O}\right)P_{\mathbb{C}}^{1}$
(\ref{(CxO)P^1(C)}) have the general structures
\begin{eqnarray}
 &  & \underset{\text{structure~}I}{\frac{\text{Spin}\left(p,10-p;\mathbb{C}\right)}{\text{Spin}\left(p,8-p;\mathbb{C}\right)\otimes U_{1}}\,}\,\text{\,\,\, or }\,\,\,\,\underset{\text{structure~}II}{\frac{\text{Spin}\left(p,10-p;\mathbb{C}\right)}{\text{Spin}\left(p-2,10-p;\mathbb{C}\right)\otimes U_{1}}},\label{real-CO-1}\\
\nonumber \\
 &  & \underset{\text{structure~}III}{\frac{\text{Spin}\left(p,10-p;\mathbb{C}\right)}{\text{Spin}\left(p-1,9-p;\mathbb{C}\right)\otimes SO\left(1,1\right)}},\label{real-CO-2}\\
\nonumber \\
 &  & Y_{16}:=\frac{\text{SO}^{\ast}\left(10\right)}{\text{SO}^{\ast}\left(8\right)\otimes\text{U}\left(1\right)}\simeq\frac{\text{SO}^{\ast}\left(10\right)}{\text{Spin}\left(6,2\right)\otimes\text{U}\left(1\right)},\label{Y-16}
\end{eqnarray}
with $p=0,1,...,10$, such that structures $I$ and $II$ get exchanged
(while structure $III$ is invariant) under $p\leftrightarrow10-p$.

The planes $\left(\mathbb{C\otimes O}\right)P^{1}$, $\left(\mathbb{C\otimes O}\right)H^{1}$,
$\left(\mathbb{C\otimes O}_{s}\right)P^{1}$ respectively correspond
to structure $I$ with $p=0$, $8$, $6$ (or structure $II$ with
$p=10$, $2$, $4$), whereas $\left(\mathbb{C}_{s}\mathbb{\otimes O}\right)P^{1}$
and $\left(\mathbb{C}_{s}\mathbb{\otimes O}_{s}\right)P^{1}$ respectively
correspond to structure $III$ with $p=9$ (or $p=1$) and $p=5$.
All other real forms of $\left(\mathbb{C}\otimes\mathbb{O}\right)P_{\mathbb{C}}^{1}$
(\ref{(CxO)P^1(C)}), and in particular (\ref{Y-16}), cannot be characterized
as Rosenfeld lines over $\mathbb{C}$(or$~\mathbb{C}_{s}$)$\otimes\mathbb{O}$(or$~\mathbb{O}_{s}$).

\subsection{$\mathbb{H}\otimes\mathbb{O}$}

In the quateroctonionic case, i.e. for $\mathbb{H}\otimes\mathbb{O}$,
the starting point is given by the \textit{quateroctonionic projective
line} over $\mathbb{C}$, i.e. by the locally symmetric coset manifold
having as isometry group $\text{Spin}\left(12,\mathbb{C}\right)$
with $\text{Spin}\left(8,\mathbb{C}\right)\otimes$Spin$\left(4,\mathbb{C}\right)$
as isotropy group : 
\begin{equation}
\left(\mathbb{H}\otimes\mathbb{O}\right)P_{\mathbb{C}}^{1}\simeq\frac{\text{Spin}\left(12,\mathbb{C}\right)}{\text{Spin}\left(8,\mathbb{C}\right)\otimes\text{Spin}\left(4,\mathbb{C}\right)}.\label{(HxO)P^1(C)}
\end{equation}
In this coset space formulation, the tangent space of $\left(\mathbb{H}\otimes\mathbb{O}\right)P_{\mathbb{C}}^{1}$
can be identified with the $\left(\boldsymbol{8}_{\left(v,s,c\right)},\boldsymbol{2},\boldsymbol{2}\right)$
representation of Spin$\left(8,\mathbb{C}\right)\otimes$Spin$\left(4,\mathbb{C}\right)$,
where we used Spin$\left(4,\mathbb{C}\right)\simeq$SL$\left(2,\mathbb{C}\right)\otimes$SL$\left(2,\mathbb{C}\right)$
and, as above, $\mathfrak{d}_{4}$-triality \cite{D4-triality} allows
to equivalently choose $v$, $s$, or $c$.

Then, by specifying the formulæ\ (\ref{eq:projectiv-line}) and (\ref{eq:hyperb-line})
for $\mathbb{A}=\mathbb{H}$ or $\mathbb{H}_{s}$ and $\mathbb{B}=\mathbb{O}$
or $\mathbb{O}_{s}$, we obtain all real forms\footnote{By real forms of $\left(\mathbb{H}\otimes\mathbb{O}\right)P_{\mathbb{C}}^{1}$,
we here mean the cosets with isometry groups given by \textit{all}
real (compact and non-compact) forms of Spin$\left(12,\mathbb{C}\right)$,
and with isotropy group given by all (compact and non-compact) real
forms of Spin$\left(8,\mathbb{C}\right)\otimes$Spin$\left(4,\mathbb{C}\right)$
which are subgroups of the corresponding real form of Spin$\left(12,\mathbb{C}\right)$.} of (\ref{(HxO)P^1(C)}) which can be expressed as Rosenfeld lines
over $\mathbb{H}$(or$~\mathbb{H}_{s}$)$\otimes\mathbb{O}$(or$~\mathbb{O}_{s}$);
they are summarized by the following Table :
\begin{center}
\begin{tabular}{|c|c|c|c|c|c|}
\hline 
\textbf{Plane} & \textbf{Isometry} & \textbf{Isotropy} & $\#_{nc}$ & $\#_{c}$ & $\chi$\tabularnewline
\hline 
\hline 
$\left(\mathbb{H}\otimes\mathbb{O}\right)P^{1}$ & $\text{Spin}\left(12\right)$ & $\text{Spin}\left(8\right)\otimes\text{SU}\left(2\right)\otimes\text{SU}\left(2\right)$ & $0$ & $32$ & $-32$\tabularnewline
\hline 
$\left(\mathbb{H}\otimes\mathbb{O}\right)H^{1}$ & $\text{Spin}\left(8,4\right)$ & $\text{Spin}\left(8\right)\otimes\text{SU}\left(2\right)\otimes\text{SU}\left(2\right)$ & $32$ & $0$ & $32$\tabularnewline
\hline 
$\left(\mathbb{H}\otimes\mathbb{O}_{s}\right)P^{1}$ & $\text{Spin}\left(8,4\right)$ & $\text{Spin}\left(4,4\right)\otimes\text{SU}\left(2\right)\otimes\text{SU}\left(2\right)$ & $16$ & $16$ & $0$\tabularnewline
\hline 
$\left(\mathbb{H}_{s}\otimes\mathbb{O}\right)P^{1}$ & $\text{Spin}\left(10,2\right)$ & $\text{Spin}\left(8\right)\otimes\text{SL}\left(2,\mathbb{R}\right)\otimes\text{SL}\left(2,\mathbb{R}\right)$ & $16$ & $16$ & $0$\tabularnewline
\hline 
$\left(\mathbb{H}_{s}\otimes\mathbb{O}_{s}\right)P^{1}$ & $\text{Spin}\left(6,6\right)$ & $\text{Spin}\left(4,4\right)\otimes\text{SL}\left(2,\mathbb{R}\right)\otimes\text{SL}\left(2,\mathbb{R}\right)$ & $16$ & $16$ & $0$\tabularnewline
\hline 
\end{tabular} 
\par\end{center}

\bigskip{}
 along with 
\begin{align}
\left(\mathbb{H\otimes O}_{s}\right)P^{1} & \simeq\left(\mathbb{H\otimes O}_{s}\right)H^{1};\\
\left(\mathbb{H}_{s}\mathbb{\otimes O}\right)P^{1} & \simeq\left(\mathbb{H}_{s}\mathbb{\otimes O}\right)H^{1};\\
\left(\mathbb{H}_{s}\mathbb{\otimes O}_{s}\right)P^{1} & \simeq\left(\mathbb{H}_{s}\mathbb{\otimes O}_{s}\right)H^{1}.
\end{align}

The real forms of $\left(\mathbb{H}\otimes\mathbb{O}\right)P_{\mathbb{C}}^{1}$
(\ref{(HxO)P^1(C)}) have the general structures
\begin{eqnarray}
 &  & \underset{\text{structure~}I}{\frac{\text{Spin}\left(p,12-p;\mathbb{C}\right)}{\text{Spin}\left(p,8-p;\mathbb{C}\right)\otimes SU(2)^{\otimes2}}}\,\,\,~\text{or}\,\,\,\,~\underset{\text{structure~}II}{\frac{\text{Spin}\left(p,12-p;\mathbb{C}\right)}{\text{Spin}\left(p-4,12-p;\mathbb{C}\right)\otimes SU(2)_{1}^{\otimes2}}},\label{real-HO-1}\\
\nonumber \\
 &  & \underset{\text{structure~}III}{\frac{\text{Spin}\left(p,12-p;\mathbb{C}\right)}{\text{Spin}\left(p-1,9-p;\mathbb{C}\right)\otimes SL(2,\mathbb{C})_{\mathbb{R}}}}~\,\,\,\text{or}~\,\,\,\,\underset{\text{structure~}IV}{\frac{\text{Spin}\left(p,12-p;\mathbb{C}\right)}{\text{Spin}\left(p-3,11-p;\mathbb{C}\right)\otimes SL(2,\mathbb{C})_{\mathbb{R}}}},\label{real-HO-2}\\
\nonumber \\
 &  & \underset{\text{structure~}V}{\frac{\text{Spin}\left(p,12-p;\mathbb{C}\right)}{\text{Spin}\left(p-2,10-p;\mathbb{C}\right)\otimes SL(2,\mathbb{R})^{\otimes2}}},\label{real-HO-3}\\
 &  & Y_{32}:=\frac{\text{SO}^{\ast}\left(12\right)}{\text{SO}^{\ast}\left(8\right)\otimes\text{SO}^{\ast}\left(4\right)}\simeq\frac{\text{SO}^{\ast}\left(12\right)}{\text{Spin}\left(6,2\right)\otimes\text{SU}\left(2\right)\otimes\text{SL}\left(2,\mathbb{R}\right)},\label{Y-32}
\end{eqnarray}
with $p=0,1,...,12$, such that structures $I$ and $II$ (and structures
$III$ and $IV$) get exchanged (while structure $V$ is invariant)
under $p\leftrightarrow12-p$.

$\left(\mathbb{H\otimes O}\right)P^{1}$, $\left(\mathbb{H\otimes O}\right)H^{1}$
and $\left(\mathbb{H\otimes O}_{s}\right)P^{1}$ respectively correspond
to structure $I$ with $p=0$, $8$, $4$ (or structure $II$ with
$p=12$, $4$, $8$), whereas $\left(\mathbb{H}_{s}\mathbb{\otimes O}\right)P^{1}$
and $\left(\mathbb{H}_{s}\mathbb{\otimes O}_{s}\right)P^{1}$ respectively
correspond to structure $V$ with $p=2$ (or $p=10$) and $p=6$.

All other real forms of $\left(\mathbb{H}\otimes\mathbb{O}\right)P_{\mathbb{C}}^{1}$
(\ref{(HxO)P^1(C)}), and in particular (\ref{Y-32}), cannot be characterized
as Rosenfeld lines over $\mathbb{H}$(or$~\mathbb{H}_{s}$)$\otimes\mathbb{O}$(or$~\mathbb{O}_{s}$).

\subsection{$\mathbb{O}\otimes\mathbb{O}$}

In the octooctonionic case, i.e. for $\mathbb{O}\otimes\mathbb{O}$,
the starting point is given by the \textit{octooctonionic projective
line} over $\mathbb{C}$, i.e. by the locally symmetric coset manifold
having as isometry group $\text{Spin}\left(16,\mathbb{C}\right)$
with $\text{Spin}\left(8,\mathbb{C}\right)\otimes$Spin$\left(8,\mathbb{C}\right)$
as isotropy group : 
\begin{equation}
\left(\mathbb{O}\otimes\mathbb{O}\right)P_{\mathbb{C}}^{1}\simeq\frac{\text{Spin}\left(12,\mathbb{C}\right)}{\text{Spin}\left(8,\mathbb{C}\right)\otimes\text{Spin}\left(8,\mathbb{C}\right)}.\label{(OxO)P^1(C)}
\end{equation}
In this coset space formulation, the tangent space of $\left(\mathbb{O}\otimes\mathbb{O}\right)P_{\mathbb{C}}^{1}$
can be identified with the $\left(\boldsymbol{8}_{\left(v,s,c\right)},\boldsymbol{8}_{\left(v,s,c\right)}\right)$
representation of Spin$\left(8,\mathbb{C}\right)\otimes$Spin$\left(8,\mathbb{C}\right)$,
where, as above, $\mathfrak{d}_{4}$-triality \cite{D4-triality}
allows to equivalently choose $v$, $s$, or $c$ (in all possible
pairs for the two Spin$\left(8,\mathbb{C}\right)$ factors of the
isotropy group).

Then, by specifying the formulæ\ (\ref{eq:projectiv-line}) and (\ref{eq:hyperb-line})
for $\mathbb{A}=\mathbb{O}$ or $\mathbb{O}_{s}$ and $\mathbb{B}=\mathbb{O}$
or $\mathbb{O}_{s}$, we obtain all real forms\footnote{By real forms of $\left(\mathbb{O}\otimes\mathbb{O}\right)P_{\mathbb{C}}^{1}$,
we here mean the cosets with isometry groups given by \textit{all}
real (compact and non-compact) forms of Spin$\left(16,\mathbb{C}\right)$,
and with isotropy group given by all (compact and non-compact) real
forms of Spin$\left(8,\mathbb{C}\right)\otimes$Spin$\left(8,\mathbb{C}\right)$
which are subgroups of the corresponding real form of Spin$\left(16,\mathbb{C}\right)$.} of (\ref{(OxO)P^1(C)}) which can be expressed as Rosenfeld lines
over $\mathbb{O}$(or$~\mathbb{O}_{s}$)$\otimes\mathbb{O}$(or$~\mathbb{O}_{s}$);
they are summarized by the following Table :
\begin{center}
\begin{tabular}{|c|c|c|c|c|c|}
\hline 
\textbf{Plane} & \textbf{Isometry} & \textbf{Isotropy} & $\#_{nc}$ & $\#_{c}$ & $\chi$\tabularnewline
\hline 
\hline 
$\left(\mathbb{O}\otimes\mathbb{O}\right)P^{1}$ & $\text{Spin}\left(16\right)$ & $\text{Spin}\left(8\right)\otimes\text{Spin}\left(8\right)$ & $0$ & $64$ & $-64$\tabularnewline
\hline 
$\left(\mathbb{O}\otimes\mathbb{O}\right)H^{1}$ & $\text{Spin}\left(8,8\right)$ & $\text{Spin}\left(8\right)\otimes\text{Spin}\left(8\right)$ & $64$ & $0$ & $64$\tabularnewline
\hline 
$\left(\mathbb{O}\otimes\mathbb{O}_{s}\right)P^{1}$ & $\text{Spin}\left(12,4\right)$ & $\text{Spin}\left(8\right)\otimes\text{Spin}\left(4,4\right)$ & $32$ & $32$ & $0$\tabularnewline
\hline 
$\left(\mathbb{O}_{s}\otimes\mathbb{O}_{s}\right)P^{1}$ & $\text{Spin}\left(8,8\right)$ & $\text{Spin}\left(4,4\right)\otimes\text{Spin}\left(4,4\right)$ & $32$ & $32$ & $0$\tabularnewline
\hline 
\end{tabular} 
\par\end{center}

\bigskip{}
 along with 
\begin{align}
\left(\mathbb{O\otimes O}_{s}\right)P^{1} & \simeq\left(\mathbb{O\otimes O}_{s}\right)H^{1},\\
\left(\mathbb{O}_{s}\mathbb{\otimes O}_{s}\right)P^{1} & \simeq\left(\mathbb{O}_{s}\mathbb{\otimes O}_{s}\right)H^{1}.
\end{align}

The real forms of $\left(\mathbb{O}\otimes\mathbb{O}\right)P_{\mathbb{C}}^{1}$
(\ref{(OxO)P^1(C)}) have the general structures
\begin{eqnarray}
 &  & \underset{\text{structure~}I}{\frac{\text{Spin}\left(p,16-p;\mathbb{C}\right)}{\text{Spin}\left(p,8-p;\mathbb{C}\right)\otimes\text{Spin}\left(8,\mathbb{C}\right)}}~\,\,\,\,\text{or}\,\,\,\,~\underset{\text{structure~}II}{\frac{\text{Spin}\left(p,16-p;\mathbb{C}\right)}{\text{Spin}\left(p-8,16-p;\mathbb{C}\right)\otimes\text{Spin}\left(8,\mathbb{C}\right)}},\label{real-OO-1}\\
\nonumber \\
 &  & \underset{\text{structure~}III}{\frac{\text{Spin}\left(p,16-p;\mathbb{C}\right)}{\text{Spin}\left(p-1,9-p;\mathbb{C}\right)\otimes\text{Spin}\left(1,7;\mathbb{C}\right)}}~\text{\,\,\,or}\,\,\,~\underset{\text{structure~}IV}{\frac{\text{Spin}\left(p,16-p;\mathbb{C}\right)}{\text{Spin}\left(p-7,15-p;\mathbb{C}\right)\otimes\text{Spin}\left(1,7;\mathbb{C}\right)}},\\
\nonumber \\
 &  & \underset{\text{structure~}V}{\frac{\text{Spin}\left(p,16-p;\mathbb{C}\right)}{\text{Spin}\left(p-2,10-p;\mathbb{C}\right)\otimes\text{Spin}\left(2,6;\mathbb{C}\right)}}~\text{\,\,\,or}~\,\,\,\underset{\text{structure~}VI}{\frac{\text{Spin}\left(p,16-p;\mathbb{C}\right)}{\text{Spin}\left(p-6,14-p;\mathbb{C}\right)\otimes\text{Spin}\left(2,6;\mathbb{C}\right)}},\\
\nonumber \\
 &  & \underset{\text{structure~}VII}{\frac{\text{Spin}\left(p,16-p;\mathbb{C}\right)}{\text{Spin}\left(p-3,11-p;\mathbb{C}\right)\otimes\text{Spin}\left(3,5;\mathbb{C}\right)}}~\text{\,\,\,or}~\underset{\text{structure~}VIII}{\,\,\,\,\frac{\text{Spin}\left(p,16-p;\mathbb{C}\right)}{\text{Spin}\left(p-5,13-p;\mathbb{C}\right)\otimes\text{Spin}\left(3,5;\mathbb{C}\right)}},\\
\nonumber \\
 &  & \underset{\text{structure~}IX}{\frac{\text{Spin}\left(p,16-p;\mathbb{C}\right)}{\text{Spin}\left(p-4,12-p;\mathbb{C}\right)\otimes\text{Spin}\left(4,4;\mathbb{C}\right)}},\label{real-OO-5}\\
\nonumber \\
 &  & Y_{64}:=\frac{\text{SO}^{\ast}\left(16\right)}{\text{SO}^{\ast}\left(8\right)\otimes\text{SO}^{\ast}\left(8\right)}\simeq\frac{\text{SO}^{\ast}\left(16\right)}{\text{Spin}\left(6,2\right)\otimes\text{Spin}\left(6,2\right)},\label{Y-64}
\end{eqnarray}
with $p=0,1,...,16$, such that structures $I$ and $II$, $III$
and $IV$, $V$ and $VI$, $VII$ and $VIII$ get exchanged (while
structure $IX$ is invariant) under $p\leftrightarrow16-p$.

The planes $\left(\mathbb{O\otimes O}\right)P^{1}$, $\left(\mathbb{O\otimes O}\right)H^{1}$
and $\left(\mathbb{O\otimes O}_{s}\right)P^{1}$ respectively correspond
to structure $I$ with $p=0$, $8$, $4$ (or structure $II$ with
$p=16$, $8$, $12$), whereas $\left(\mathbb{O}_{s}\mathbb{\otimes O}_{s}\right)P^{1}$
corresponds to structure $IX$ with $p=8$. All other real forms of
$\left(\mathbb{O}\otimes\mathbb{O}\right)P_{\mathbb{C}}^{1}$ (\ref{(OxO)P^1(C)}),
and in particular (\ref{Y-64}), cannot be characterized as Rosenfeld
lines over $\mathbb{O}$(or$~\mathbb{O}_{s}$)$\otimes\mathbb{O}$(or$~\mathbb{O}_{s}$).

\section{\label{5}Conclusions}

In this work, in order to provide a rigorous characterization of Rosenfeld
``projective''\ spaces over (rank-2) tensor products of (division
or split) Hurwitz algebras, we have introduced some ``magic''\ formulæ\ (\ref{eq:projectiv})-(\ref{eq:mixed})
(for planes) and (\ref{eq:projectiv-line})-(\ref{eq:hyperb-line})
(for lines), which allowed us to characterize many Rosenfeld spaces
as symmetric (pseudo-)Riemannian cosets, whose isometry and isotropy
groups have Lie algebras which are entries of real forms of the order-3
(Freudenthal-Tits \cite{Tits,Freud-1965}) Magic Square, or of the
order-2 (Barton-Sudbery \cite{BS}) Magic Square.

For ``projective''\ planes, the application of the ``magic''\ formulæ\ (\ref{eq:projectiv})-(\ref{eq:mixed})
to the case in which at least one of the two Hurwitz algebras in the
associated tensor product is given by the octonions $\mathbb{O}$
or by the split octonions $\mathbb{O}_{s}$, allows us to retrieve
all (compact and non-compact) real forms of the corresponding \textit{octonionic}
Rosenfeld planes, \textit{except} for a limited numbers of pseudo-Riemannian
symmetric cosets, named $X_{32,I}$, $X_{32,II}$, $X_{64,I}$, $X_{64,II}$
,$X_{64,III}$, $X_{128,I}$ and $X_{128,II}$, and respectively given
by (\ref{32-1}), (\ref{32-2}), (\ref{64-1}), (\ref{64-2}), (\ref{64-3}),
(\ref{128-1}) and (\ref{128-2}). All such cosets share a common
property : up to some possible rank-1 Lie group factor, their isotropy
group is given by the non-compact Lie group SO$^{\ast}\left(N\right)$
for $N=10,12,16$. The fact that such spaces are not encompassed by
our ``magic''\ formulæ\ (\ref{eq:projectiv})-(\ref{eq:mixed})
is ultimately due to the fact that the corresponding Lie algebra $\mathfrak{so}^{\ast}(N)$
does not occur in any real form of the Magic Square of order 2 (except
for the case $N=8$, for which it holds the special isomorphism $\mathfrak{so}^{\ast}(8)\simeq\mathfrak{so}(6,2)$).

For ``projective''\ lines, the application of the ``magic''\ formulæ\ (\ref{eq:projectiv-line})-(\ref{eq:hyperb-line})
to the case in which at least one of the two Hurwitz algebras in the
associated tensor product is given by the octonions $\mathbb{O}$
or by the split octonions $\mathbb{O}_{s}$, allows us to retrieve
some (compact and non-compact) real forms of the corresponding \textit{octonionic}
Rosenfeld lines, \textit{but still a number of other real forms of
Rosenfeld lines is left out}. Again, there are spaces, such as $Y_{16},Y_{32}$
and $Y_{64}$ (respectively given by (\ref{Y-16}), (\ref{Y-32})
and (\ref{Y-64})) that have their isometry and isotropy groups containing
factors SO$^{\ast}(M)$ for $M=4,8,10,12,16$; but there are also
other pseudo-Riemannian spaces left out, given by (\ref{real-O}),
(\ref{real-CO-1})-(\ref{real-CO-2}), (\ref{real-HO-1})-(\ref{real-HO-3})
and (\ref{real-OO-1})-(\ref{real-OO-5}).

The fact that some pseudo-Riemannian symmetric spaces are not covered
by the classification yielded by the ``magic'' formulæ (\ref{eq:projectiv})-(\ref{eq:mixed})
and (\ref{eq:projectiv-line})-(\ref{eq:hyperb-line}), which are
related to the entries of the Magic Squares of order 2 and 3, means
that \textit{not all} real forms of Rosenfeld spaces can be realized
as ``projective''\ lines or planes over (rank 2) tensor products
of (division or split) Hurwitz algebras. If one still believe that
Rosenfeld's approach was right after all, one might want to extend
our ``magic''\ formulæ\ to involve not only \textit{unital, alternative}
composition\textit{\ }algebras, i.e. (division or split) Hurwitz
algebras, but also other algebras which are not unital or alternative.
Extensive work on Magic Squares over flexible composition algebras
has been done by Elduque (see e.g. \cite{EldMS1,EldMS2}), whereas
a geometrical framework for these algebras has been recently discussed
in \cite{CorrZucc} and \cite{CorrMarZucc}.

\textit{Therefore, after all, Rosenfeld might still be right}....we
hope to report on this in forthcoming investigations.

\section*{Acknowledgments}

The work of D. Corradetti is supported by a grant of the Quantum Gravity
Research Institute. The work of AM is supported by a ``Maria Zambrano''\ distinguished
researcher fellowship at the University of Murcia, ES, financed by
the European Union within the NextGenerationEU program.

\newpage{}

\bigskip{}

\end{document}